\pdfoutput=1
\RequirePackage{ifpdf}
\ifpdf 
\documentclass[pdftex]{sigma}
\else
\documentclass{sigma}
\fi

\usepackage{tikz}
\usepackage{mathrsfs}
\usepackage{mathtools,leftidx}
\usepackage[shortcuts]{extdash}

\numberwithin{equation}{section}

\newtheorem{Theorem}{Theorem}[section]
\newtheorem{Lemma}[Theorem]{Lemma}

\newtheorem{Corollary}[Theorem]{Corollary}
\newtheorem{Proposition}[Theorem]{Proposition}
\newtheorem{Claim}[Theorem]{Claim}

 { \theoremstyle{definition}
\newtheorem{Definition}[Theorem]{Definition}
\newtheorem{Remark}[Theorem]{Remark} }

\usepackage{gasymbols}

\makeatletter

\newcounter{case}
\@addtoreset{case}{proof}

\makeatother

\makeatletter
\newcommand{\leqnmode}{\tagsleft@true}
\makeatother
\renewcommand{\ncapa}{\mathfrak{c}}

\newcommand{\eum}{\delta}
\newcommand{\ma}{\mathfrak{m}}

\newcommand{\tp}[1][p]{{\text{\tiny$(#1)$}}}

\newcommand{\lop}{\mathscr{L}}

\begin{document}

\allowdisplaybreaks

\renewcommand{\thefootnote}{}

\newcommand{\arXivNumber}{2305.01453}

\renewcommand{\PaperNumber}{091}

\FirstPageHeading

\ShortArticleName{Nonlinear Isocapacitary Concepts of Mass}

\ArticleName{Nonlinear Isocapacitary Concepts of Mass \\ in 3-Manifolds with Nonnegative Scalar Curvature\footnote{This paper is a~contribution to the Special Issue on Differential Geometry Inspired by Mathematical Physics in honor of Jean-Pierre Bourguignon for his 75th birthday. The~full collection is available at \href{https://www.emis.de/journals/SIGMA/Bourguignon.html}{https://www.emis.de/journals/SIGMA/Bourguignon.html}}}

\Author{Luca BENATTI~$^{\rm a}$, Mattia FOGAGNOLO~$^{\rm b}$ and Lorenzo MAZZIERI~$^{\rm c}$}

\AuthorNameForHeading{L.~Benatti, M.~Fogagnolo and L.~Mazzieri}

\Address{$^{\rm a)}$~Universit\`a degli Studi di Pisa, Largo Bruno Pontecorvo 5, 56127 Pisa, Italy}
\EmailD{\href{mailto:luca.benatti@unipi.it}{luca.benatti@unipi.it}}
\URLaddressD{\url{https://sites.google.com/view/luca-benatti}}

\Address{$^{\rm b)}$~Universit\`a di Padova, via Trieste 63, 35121 Padova, Italy}
\EmailD{\href{mailto:mattia.fogagnolo@unipd.it}{mattia.fogagnolo@unipd.it}}
\URLaddressD{\url{https://sites.google.com/view/mattiafogagnolo}}

\Address{$^{\rm c)}$~Universit\`a degli Studi di Trento, via Sommarive 14, 38123 Povo (TN), Italy}
\EmailD{\href{mailto:lorenzo.mazzieri@unitn.it}{lorenzo.mazzieri@unitn.it}}
\URLaddressD{\url{https://sites.google.com/site/mazzierihome}}

\ArticleDates{Received May 03, 2023, in final form October 23, 2023; Published online November 10, 2023}

\Abstract{We deal with suitable nonlinear versions of Jauregui's isocapacitary mass in $3$-manifolds with nonnegative scalar curvature and compact outermost minimal boundary. These masses, which depend on a parameter $1<p\leq 2$, interpolate between Jauregui's mass~${p=2}$ and Huisken's isoperimetric mass, as $p \to 1^+$. We derive positive mass theorems for these masses under mild conditions at infinity, and we show that these masses do coincide with the $\ADM$ mass when the latter is defined. We finally work out a nonlinear potential theoretic proof of the Penrose inequality in the optimal asymptotic regime.}


\Classification{83C99; 35B40; 35A16; 31C15; 53C21}

\begin{flushright}
\begin{minipage}{65mm}
\em Dedicated to Jean–Pierre Bourguignon\\
on the occasion of his 75th birthday
\end{minipage}
\end{flushright}

\renewcommand{\thefootnote}{\arabic{footnote}}
\setcounter{footnote}{0}

\section{Introduction}
The isoperimetric concept of mass was introduced by Huisken~\cite{huisken_isoperimetricconceptmassgeneral_2006} to study $3$-manifolds with nonnegative scalar curvature that are possibly nonsmooth, and not fulfilling the asymptotic requirements needed to define the classic $\ADM$ mass~\cite{arnowitt_coordinateinvarianceenergyexpressions_1961}. Given a $3$-manifold $(M,g)$, the isoperimetric mass is indeed defined as
\begin{equation}\label{eq:isoperimetric_mass_intro}
 \ma_{\iso}=\sup_{(\Omega_j)_{j\in \N}}\limsup_{j \to +\infty} \ma_{\iso}(\Omega_j),
\end{equation}
where the supremum is taken among all exhaustions $(\Omega_j)_{j\in \N}$ consisting of domains with $\CS^{1}$-boundary and the \textit{isoperimetric quasi-local mass} $\ma_{\iso}(\Omega)$ is
\begin{equation*}
 \ma_{\iso}(\Omega)=\frac{2}{\abs{\partial \Omega}}\bigg( \abs{\Omega} -\frac{\abs{\partial \Omega}^{\frac{3}{2}}}{6 \sqrt{\pi}}\bigg).
\end{equation*}
Notice that, without any further assumption on the $3$-manifold, the isoperimetric mass might in principle be any number in $[-\infty,+\infty]$. Unlike the $\ADM$ mass, which is defined on an asymptotically flat chart $x=\big(x^1,x^2,x^3\big)$ at infinity as
\begin{equation*}
 \ma_{\ADM} = \lim_{r \to +\infty}\frac{1}{16\pi} \int_{\set{\abs{x}= r}}(\partial_k g_{ii}- \partial_i g_{ki} ) \frac{x^k}{\abs{x}}\dif \sigma,
\end{equation*}
the isoperimetric mass does not require passing to a chart to be defined. Rather, it is based on the geometric concepts of volume and perimeter, making it well-defined even when there is limited information on the asymptotic behaviour of the metric.
Inspired by this observation, in~\cite{benatti_isoperimetricriemannianpenroseinequality_2022}, we proved a Riemannian Penrose inequality~\cite{bray_proofriemannianpenroseinequality_2001, huisken_inversemeancurvatureflow_2001} for the isoperimetric mass in the class of \emph{strongly $1$-nonparabolic} Riemannian manifolds with nonnegative scalar curvature. With the locution {\em strongly $1$-nonparabolic} manifolds, we denote manifolds on which any bounded $\Omega \subset M$, whose boundary is homologous to $\partial M$, admits a proper locally Lipschitz weak inverse mean curvature flow (IMCF for short), that is a solution $w_1$ to the problem
\begin{equation}\label{eq:IMCF}
 \begin{cases}
 \div\bigg( \dfrac{ \D w_1}{\abs{\D w_1}}\bigg) = \abs{ \D w_1} & \text{on $M\smallsetminus \Omega$,}\\
 w_1 = 0 & \text{on $\partial \Omega$,}\\
 w_1 \to +\infty & \text{as $\dist(x,\partial \Omega) \to +\infty$},
 \end{cases}
\end{equation}
according to the definition introduced in~\cite{huisken_inversemeancurvatureflow_2001}. The analysis leading to the isoperimetric Riemannian Penrose inequality in~\cite{benatti_isoperimetricriemannianpenroseinequality_2022} was carried out using a new asymptotic comparison between the Hawking mass (see~\eqref{eq:hawking-mass} below) and the isoperimetric mass along the level sets of the weak IMCF.

In the present paper, we are going to develop a similar theory, in the case where the weak IMCF~\eqref{eq:IMCF} is replaced by the level set flow of weak solutions $w_p\in \CS^{1,\beta}_{\loc}(M \smallsetminus \Omega)$ to the boundary value problem
\begin{equation}
\label{eq:pb-intro}
 \begin{cases}
 \Delta_p w_p = \abs{\D w_p}^p &\text{on $M \smallsetminus \Int \Omega$,}\\
 w_p=0 & \text{on $\partial \Omega$,}\\
 w_p\to +\infty & \text{as $\dist(x,\partial \Omega) \to +\infty$.}
 \end{cases}
\end{equation}
The link between the above problem and the weak IMCF~\eqref{eq:IMCF} relies on the fact that $w_p\to w_1$ as $p\to 1^+$ locally uniformly on $M$~\cite{kotschwar_localgradientestimatesharmonic_2009,mari_flowlaplaceapproximationnew_2022, moser_inversemeancurvatureflow_2007,moser_inversemeancurvatureflow_2008}, provided some natural global requirements are met by the manifold $(M, g)$. On the other hand,
the solutions to problem~\eqref{eq:pb-intro} are directly connected to the notion of $p$-capacitary potential of a compact body $\Omega$. In fact, setting $w_p = -(p-1) \log u_p$ implies that $u_p$ is $p$-harmonic, that is $\Delta_p u_p = 0$. These relationships have been instrumental in demonstrating a series of geometric inequalities by means of monotonicity formulas, holding along the level sets of solutions to equation~\eqref{eq:pb-intro}. As $p \to 1^+$, these inequalities become increasingly close to the desired result.
This machinery, firstly introduced in the case $p=2$ for harmonic functions in~\cite{agostiniani_sharpgeometricinequalitiesclosed_2020, agostiniani_monotonicityformulaspotentialtheory_2020}, has proven to be powerful enough to produce an enhanced version of the Minkowski inequality~\cite{agostiniani_minkowskiinequalitiesnonlinearpotential_2022, fogagnolo_geometricaspectscapacitarypotentials_2019}, later extended to Riemannian manifolds with nonegative Ricci curvature~\cite{benatti_minkowskiinequalitycompleteriemannian_2022} as well as to the anisotropic setting~\cite{xia_anisotropiccapacityanisotropicminkowski_2022}.
 In~\cite{agostiniani_greenfunctionproofpositive_2021} and subsequently in~\cite{agostiniani_riemannianpenroseinequalitynonlinear_2022}, the authors used this approach on $3$-manifolds with nonnegative scalar curvature to prove the Riemannian Penrose inequality for a single black hole, based on the monotonic behaviour of a suitable $p$-harmonic version of the Hawking mass (see~\eqref{eq:p-Hawkingmass} below).
\emph{Throughout the manuscript, Riemannian manifolds are assumed to be smooth, connected, metrically complete, noncompact, with one single end.}

The main object of interest in the present paper is the following nonlinear potential theoretic version of Huisken's isoperimetric mass~\eqref{eq:isoperimetric_mass_intro}, that we call the {\em $p$-isocapacitary mass}. It involves the classical notion of $p$-capacity of a compact set $K \subset M$, that, in dimension $3$ and for $1 < p <3$, is given by
\begin{equation*}
 \ncapa_p(K) = \inf \set{\frac{1}{4\pi}\bigg(\frac{p-1}{3-p}\bigg)^{p-1}\int_{M\smallsetminus K} \abs{\D v}^p \dif \mu\st v \in \CS_c^\infty(M),\, v \geq 1\text{ on } K}.
\end{equation*}
When the boundary of $K$ is regular enough, such infimum is realized by the $p$-capacitary potential of $K$, i.e., the unique $p$-harmonic function $u_p$, equal to $1$ on $\partial K$ and vanishing at infinity.
\begin{Definition}[$p$-isocapacitary mass]
Let $(M,g)$ be a Riemannian $3$-manifold with compact boundary, and let $1<p<3$. Given a closed bounded subset $\Omega\subset M$ containing $\partial M$ with $\CS^{1}$-boundary the quasi-local $p$-isocapacitary mass of $\Omega$ is defined as
\begin{equation*}
 \ma_{\iso}^{\tp}(\Omega)=\frac{1}{2p \pi \ncapa_p(\partial \Omega)^{\frac{2}{3-p}}}\left( \abs{\Omega} -\frac{4 \pi}{3} \ncapa_p (\partial \Omega)^\frac{3}{3-p}\right).
\end{equation*}
The $p$-isocapacitary mass $\ma^{\tp}_{\iso}$ of $(M,g)$ is defined as \begin{equation}\label{eq:isopcapacitary_mass}
 \ma_{\iso}^{\tp} = \sup_{(\Omega_j)_{j \in \N}} \limsup_{j \to +\infty} \ma_{\iso}^{\tp}(\partial \Omega_j),
\end{equation}
where the supremum is taken among all exhaustions $\set{\Omega_j}_{j\in \N}$.
\end{Definition}
As for the isoperimetric mass, the $p$-isocapacitary mass might be any number in $[-\infty,+\infty]$. The special and particularly relevant case of the $2$-isocapacitary mass has been recently introduced and studied by Jauregui~\cite{jauregui_admmasscapacityvolumedeficit_2020}.

A first natural and fundamental question about these newly introduced quantities is whether they are nonnegative, on the class of $3$-manifolds with nonnegative scalar curvature, where a~solution to~\eqref{eq:pb-intro} exists for any bounded $\Omega$ with regular boundary. The latter mentioned property will be referred to as the \textit{strong $p$-nonparabolicity} of the manifold $(M,g)$, a terminology that interpolates between the notion of strong nonparabolicity introduced by Ni~\cite{ni_meanvaluetheoremsmanifolds_2007} and the notion of strong $1$-nonparabolicity, which was employed in~\cite{benatti_isoperimetricriemannianpenroseinequality_2022} and recalled above.

Our first main result is a nonlinear potential-theoretic version of the Riemannian Penrose inequality, that, although not sharp, implies the positive mass theorem for the $p$-isocapacitary mass, with its associated rigidity statement. We prove its validity under the following asymptotic integral gradient estimate:
\bgroup\leqnmode
\begin{equation*}
\begin{minipage}[t]{.9585\textwidth}%
Given any $\Omega \subset M$ closed bounded subset with smooth and connected boundary homologous to $\partial M$, the function $w_p\in \CS^1_{\loc}\big(M \smallsetminus \overline{\Omega}\big)$ that solves~\eqref{eq:pb-intro} satisfies
\begin{gather*}
\int_{\partial \Omega_t} \abs{ \D w_p}^2 \dif \sigma = o\bigl(\ee^{t/(p-1)}\bigr)
\end{gather*}
as $t \to +\infty$ where $\Omega_t= \set{ w_p \leq t}$.
\end{minipage}\tag{$\dag$}\label{eq:energy_hp}
\end{equation*}
\egroup

\begin{Theorem}[$p$-isocapacitary Riemannian Penrose inequality]
\label{thm:pmt-intro}
 Let $(M,g)$ be a strongly $p$-nonparabolic Riemannian $3$-manifold, $1<p<3$, with nonnegative scalar curvature and with smooth, compact, connected, minimal, possibly empty boundary. Assume also that $(M,g)$ satisfies~\eqref{eq:energy_hp} and $H_2( M, \partial M; \Z)=\set{0}$. Then
 \begin{equation}
\label{eq:p-penroseintro-crociata}
\ncapa_p(\partial M)^{\frac{1}{3-p}}\leq \frac{5-p}{2}\ma^{\tp}_{\iso}.
\end{equation}
In particular, $\ma^{\tp}_{\iso}\geq 0$ and it vanishes if and only if $(M,g)$ is isometric to the flat $3$-dimensional Euclidean space.
\end{Theorem}
As we are going to show below in Theorem~\ref{thm:equivalence-masses}, under natural assumptions of asymptotic flatness, the $p$-capacitary masses coincide with one another, and they are all equal to the classical $\ADM$ mass. Then, in the limit as $p \to 1^+$,~\eqref{eq:p-penroseintro-crociata} yields the Riemannian Penrose inequality.
In~\cite{bray_capacitysurfacesmanifoldsnonnegative_2008} and its extension~\cite{xiao_harmoniccapacityasymptoticallyflat_2016}, the version of~\eqref{eq:p-penroseintro-crociata} with the sharp exponent is derived employing the weak inverse mean curvature flow. The above result on the other hand just involves nonlinear potential theoretic concepts, such as the notions of $p$-nonparabolicity and of $p$-isocapacitary mass.

The condition~\eqref{eq:energy_hp} is actually very mild. As we are going to detail in Remark~\ref{rmk:geometric_energy}, it is always fulfilled on manifolds that are merely $\CS^0$-asymptotically flat, provided a suitable Ricci lower bound is also satisfied. $\CS^1$-asymptotically flat Riemannian manifolds are also fulfilling condition~\eqref{eq:energy_hp} for $1<p\leq 2$, as proved in Lemma~\ref{prop:gradientdecay}. This latter class of manifolds is particularly natural in the framework of mathematical general relativity, as the works of Bartnik~\cite{bartnik_massasymptoticallyflatmanifold_1986} and Chru\'sciel~\cite{chrusciel_boundaryconditionsspatialinfinity_1986} showed that the $\ADM$ mass is well defined on $\CS^1_\tau$-asymptotically flat Riemannian manifolds, with $\tau >1/2$. In fact, our second main result shows that on $\CS^1_\tau$-asymptotically Flat Riemannian $3$-manifolds with nonnegative scalar curvature the $p$-isocapacitary masses do coincide with the $\ADM$ mass for any $1 \leq p \leq 2$. This fact was previously known only for $p = 2$ and only for harmonically flat manifolds, as proven by Jauregui in the insightful paper~\cite[Corollary~8]{jauregui_admmasscapacityvolumedeficit_2020}. Most of the authors assume that the scalar curvature belongs to $L^1(M)$ in the definition of the $\ADM$ mass so that the latter is not only a geometric invariant but also a finite number. Here, we do not assume any integrability of the scalar curvature, hence, our $\ADM$ mass is not finite a priori. The statement below, can be understood in the sense that all masses are infinite as soon as one of them is, and they all coincide with one another as soon as one of them is finite.

\begin{Theorem}\label{thm:equivalence-masses}
 Let $(M,g)$ be a $\CS^1_\tau$-asymptotically flat Riemannian $3$-manifold, $\tau >1/2$, with nonnegative scalar curvature and with smooth, compact, minimal, possibly empty boundary. Assume that $H_2(M,\partial M ; \Z)=\set{0}$. Then
 \[
 \ma^{\tp}_{\iso}=\ma_{\iso} = \ma_{\ADM}
 \]
 for all $1<p\leq 2 $.
\end{Theorem}
What we will actually prove in this paper is the chain of inequalities $\ma_{\ADM} \leq \ma^{\tp}_{\iso} \leq \ma_{\iso}$. This will in turn give the identities above due to the equality $\ma_{\iso} = \ma_{\ADM}$, obtained in the setting of Theorem~\ref{thm:equivalence-masses} in~\cite[Theorem 4.13]{benatti_isoperimetricriemannianpenroseinequality_2022}.
In the proof of Theorem~\ref{thm:equivalence-masses}, the inequality ${\ma^{\tp}_{\iso}\geq \ma_{\iso}}$ is deduced substantially arguing as in~\cite[Theorem 5]{jauregui_admmasscapacityvolumedeficit_2020}, combining some of the computations in~\cite{fan_largespheresmallspherelimitsbrownyork_2009} together with an extension of the main estimate in~\cite{bray_capacitysurfacesmanifoldsnonnegative_2008,xiao_harmoniccapacityasymptoticallyflat_2016}
for the isoperimetric mass (see Proposition~\ref{thm:p-braymiao}).

The reverse inequality, for which the harmonically flat condition was invoked in~\cite{jauregui_admmasscapacityvolumedeficit_2020}, is instead proven by integrating the sharp isoperimetric inequality in~\cite[Corollary~C.3]{chodosh_isoperimetryscalarcurvaturemass_2021} to obtain a sharp $p$-isocapacitary inequality in terms of the isoperimetric mass, Theorem~\ref{thm:sharp-isopcapacitary}. This last step is inspired by the classical derivation of the sharp $p$-isocapacitary inequality from the sharp isoperimetric inequality, as in~\cite[Theorem 4.1]{benatti_minkowskiinequalitycompleteriemannian_2022}.
The identification with the $\ADM$ mass finally follows from~\cite[Theorem 4.13]{benatti_isoperimetricriemannianpenroseinequality_2022}, where it was showed to coincide with the isoperimetric mass in the above optimal regime, sharpening~\cite[Theorem 3]{jauregui_lowersemicontinuitymass0_2019}.

It is natural to conjecture that the equivalence among $p$-isocapacitary masses also holds under weaker asymptotic assumptions, where the $\ADM$ mass may not even be well defined.
In this direction, in the generality of $\CS^0$-asymptotically flat Riemannian manifold satisfying~\eqref{eq:energy_hp} for $1<p\leq 2$, we will prove the following two-sided estimate (see Lemma~\ref{cor:lowerbound_miso} and Theorem~\ref{thm:upperbound_miso}):
\[
\ma^{\tp}_{\iso} \leq \ma_{\iso} \leq \left( \frac{2^{2p-1}\pi^{\frac{p-1}{2}} p^{p}\kst_{S}^{\frac{3}{2}(p-1)}}{(3-p) (p-1)^{p-1}} \right)^{\frac{1}{3-p}} \ma^{\tp}_{\iso},
\]
where $\kst_{S}$ is the global Sobolev constant on $(M,g)$. The lower bound is optimal, while the upper bound sharpens as $p\to 1^+$. As a consequence, in this generality $\ma_{\iso}$ can at least be recovered as the limit of its $p$-capacitary versions as $p \to 1^+$.

In concluding the paper, we propose an alternative proof of the Riemannian Penrose inequality in the sharp asymptotic regime given in~\cite{benatti_isoperimetricriemannianpenroseinequality_2022}, that is for $\CS^1_\tau$-asymptotically flat Riemannian $3$-manifolds, with $\tau > 1/2$. In this previous work, we exploited the better asymptotic behaviour of harmonic functions and the monotonicity of the $2$-Hawking mass discovered in~\cite{agostiniani_greenfunctionproofpositive_2021} to improve the original argument by Huisken and Ilmanen~\cite{huisken_inversemeancurvatureflow_2001} based on the IMCF, as far as the asymptotic analysis at infinity is concerned.

Replacing the IMCF with the level sets flow of the solutions $w_p$ to~\eqref{eq:pb-intro}, we obtain a nonsharp family of $p$-Penrose Inequalities in terms of the $p$-capacity of the horizon. These then provide the optimal and classical Riemannian Penrose inequality in the limit as $p \to 1^+$. We recall that a minimal boundary $\partial M$ is outermost if no other closed minimal surface homologous to~$\partial M$ is contained in $M$. Moreover, we understand that, in case such boundary is empty, then no minimal closed hypersurfaces exists in $(M, g)$. We can now state the last main result of the paper.
\begin{Theorem}
\label{thm:penroseadmintro}
 Let $(M,g)$ be a $\CS^1_\tau$-asymptotically flat Riemannian $3$-manifold, $\tau >1/2$, with nonnegative scalar curvature and with smooth, compact, minimal, connected and outermost possibly empty boundary. Then
\begin{equation}
\label{eq:p-penroseintro}
 \ncapa_p(\partial M)^{\frac{1}{3-p}}\leq 2 \ma_{\ADM}
\end{equation}
for any $1< p \leq 2$. Letting $p\to 1^+$, we get
\begin{equation}
\label{eq:penroseintro}
 \sqrt{\frac{\abs{ \partial M}}{16 \pi}} \leq \ma_{\ADM}.
\end{equation}
\end{Theorem}
Both in the above result and in Theorem~\ref{thm:equivalence-masses}, the restriction to $p\leq 2$ is due to technical reasons, that can be devised in the proof of Theorem~\ref{thm:upperbound_miso}.
\subsection*{Outline of the paper}In Section~\ref{sec:2}, we gather some basic facts about $p$-harmonic potentials. The content of this section is substantially well known. In~Section~\ref{sec:3}, we work out the main asymptotic comparison at infinity between the $p$-Hawking mass and the quasi-local $p$-isocapacitary mass, see Lemma~\ref{thm:mass_control_pcap}. We deduce the nonsharp Riemannian Penrose inequality Theorem~\ref{thm:pmt-intro} for the $p$-isocapacitary mass. In~Section~\ref{sec:4}, we show relations among the $p$-isocapacitary masses for the various values of $p$, in turn obtaining Theorem~\ref{thm:equivalence-masses}. Finally, in~Appendix~\ref{app:monotonicities}, we include a proof of the full monotonicity of the $p$-Hawking mass, since the original~\cite[Theorem~1.1]{agostiniani_riemannianpenroseinequalitynonlinear_2022} actually yields such result only along regular values. We will also relate such quantity with a similar one considered in~\cite{chan_monotonicitygreenfunctions_2022} that will naturally appear in the asymptotic comparison argument ruling our main results.

\section{Preliminaries in nonlinear potential theory}
\label{sec:2}
We first remind that the operator $\Delta_p$ is defined by $\Delta_p f = \div \big(\abs{\nabla f}^{p-2}\nabla f\big)$, for $f \in \CS^2$ with nonvanishing gradient.
As far as basic principles and regularity for \emph{weakly} $p$-harmonic functions are concerned, we just refer the reader to
\cite{dibenedetto_alphalocalregularityweak_1983,evans_newprooflocal1_1982,heinonen_asuperharmonicfunctionssupersolutionsdegenerate_1988,ladyzenskaja_linearquasilinearellipticequations_1968,lewis_regularityderivativessolutionscertain_1983,lieberman_boundaryregularitysolutionsdegenerate_1988,serrin_localbehaviorsolutionsquasilinear_1964,tolksdorf_dirichletproblemquasilinearequations_1983,trudinger_harnacktypeinequalitiestheir_1967,uraltseva_degeneratequasilinearellipticsystems_1968,valtorta_laplaceoperatorriemannianmanifolds_2013}, and~\cite{holopainen_nonlinearpotentialtheoryquasiregular_1990,holopainen_volumegrowthgreenfunctions_1999} for the theory of $p$-capacitary potentials in exterior domains (see also~\cite[Chapter 1]{benatti_monotonicityformulasnonlinearpotential_2022} and reference therein), including existence issues. The function $w_p$ that solves~\eqref{eq:pb-intro} can be written as $w_p= -(p-1)\log u_p$, where $u_p \in \CS^{1,\beta}_{\loc}(M \smallsetminus \Int \Omega)$ is the solution to
\begin{equation}\label{eq:p-cap-potential}
 \begin{cases}
 \Delta^{\tp}_g u_p = 0 &\text{on $M \smallsetminus \Int \Omega$,}\\
 u_p=1 & \text{on $\partial \Omega$,}\\
 u_p\to 0 & \text{as $\dist(x,\partial \Omega)\to +\infty$,}
 \end{cases}
\end{equation}
the following definition of strongly $p$-nonparabolic Riemannian manifold is consistent with that of strong nonparabolicity~\cite{ni_meanvaluetheoremsmanifolds_2007} and with the limit case of strong $1$-nonparabolicity~\cite{benatti_isoperimetricriemannianpenroseinequality_2022}.

\begin{Definition}[strongly $p$-nonparabolic]
We say that a $3$-dimensional Riemannian manifold~$(M,g)$ with compact possibly empty boundary is \emph{strongly $p$-nonparabolic}, $1<p<3$, if there exists a solution to~\eqref{eq:pb-intro} for some $\Omega \subseteq M$ closed bounded with smooth boundary homologous to $\partial M$.
\end{Definition}

\begin{Remark}
By the maximum principle, in a strongly $p$-nonparabolic manifold every $\Omega$ with $\CS^{1}$-boundary homologous to $\partial M$ admits a solution to~\eqref{eq:pb-intro}.
\end{Remark}

This definition naturally comes with the notion of the $p$-capacity of a compact subset $K \subset M$, which we recall to be
\begin{equation*}
 \ncapa_p(K) = \inf \set{\frac{1}{4\pi}\bigg(\frac{p-1}{3-p}\bigg)^{p-1}\int_{M\smallsetminus K} \abs{\D v}^p \dif \mu\st v \in \CS_c^\infty(M),\, v \geq 1\text{ on } K}.
\end{equation*}
When the boundary of $K$ is sufficiently regular, the above infimum is realized by the function~$u_p$ solving~\eqref{eq:p-cap-potential}. A fundamental property of the $p$-capacity is that it is monotone with respect to the standard inclusion of sets. More specifically, if we consider sublevel sets $\Omega_t=\set{w_p \leq t}$ of solutions to~\eqref{eq:pb-intro}, we have that their $p$-capacities grow exponentially with respect to the arrival time parameter $t$. This is
completely analogous to the exponential growth of the area along the IMCF. We recall this useful property in the following lemma, proved in~\cite[Lemma~3.8]{holopainen_nonlinearpotentialtheoryquasiregular_1990}.

\begin{Lemma}
Let $(M,g)$ be a $3$-dimensional Riemannian manifold with compact possibly empty boundary $\partial M$. Let $\Omega\subseteq M$ be a closed bounded subset homologous to $\partial M$ with $\CS^1$-boundary, $w_p$ the solution to~\eqref{eq:pb-intro} starting at $\Omega$ and $\Omega_t= \set{w_p \leq t}$, $1<p<3$. We have
\[
 \ncapa_p(\partial \Omega_t) = \ee^{t} \ncapa_p(\partial \Omega)= \frac{1}{4\pi}\int_{\partial \Omega_t} \bigg(\frac{ \abs{ \D w_p}}{3-p}\bigg)^{p-1}\dif \sigma.
\]
\end{Lemma}

\subsection{Estimates on asymptotically flat Riemannian manifolds}
In~\cite{benatti_asymptoticbehaviourcapacitarypotentials_2022}, we prove the asymptotic behaviour of $p$-harmonic potential~\eqref{eq:p-cap-potential} assuming a lower bound on the Ricci curvature other than asymptotic assumptions on the metric. Here we want to remove the additional assumption on Ricci curvature. We will prove that it is superfluous if we assume the metric is $\CS^1$-asymptotically flat.

We start by giving the precise definition of asymptotically flat $3$-manifolds.
\begin{Definition}[asymptotically flat Riemannian manifolds] 
A $3$-dimensional Riemannian manifold $(M,g)$ with compact possibly empty boundary is \emph{$\CS^{k}_\tau$-asymptotically flat} \big(resp. \emph{$\CS^{k}$-asymp\-tot\-ically flat}\big), with order $k\in \N$ and rate $\tau >0$, if the following conditions are satisfied.
\begin{enumerate}\itemsep=0pt
 \item There exists a compact set $K \subseteq M$ such that $M \smallsetminus K$ is differmorphic to $\R^3\smallsetminus \set{\abs{x}\leq R}$, through a map $\big(x^1,x^2,x^3\big)$ whose components are called \emph{asymptotically flat coordinates}.
 \item In the chart $\big(M \smallsetminus K, \big(x^1,x^2,x^3\big)\big)$ the metric tensor is expressed as
 \[
 g= g_{ij} \dd x^i \otimes\dd x^j= (\eum_{ij}+\eta_{ij}) \dd x^i \otimes\dd x^j
 \]
 with
 \[
 \sum_{i,j=1}^{3}\sum_{\abs{\beta}=0}^k \abs{x}^{\abs{\beta}+\tau} \abs{ \partial_\beta \eta_{ij}} =O(1) \quad \text{(resp. $=o(1)$), \ \ as $\abs{x} \to +\infty$.}
 \]
\end{enumerate}
\end{Definition}

$\CS^0$-asymptotically flat Riemannian manifolds are in particular strongly $p$-nonparabolic. This is a consequence of~\cite[Theorem 3.6]{mari_flowlaplaceapproximationnew_2022}, implying it under the mere existence of a global (even weighted) Sobolev inequality.

We first point out a decay estimate for the gradient of $u_p$ holding on $\CS^1$-asymptotically flat Riemannian manifolds. It is immediately obtained as a consequence of the Cheng--Yau inequality for $p$-harmonic functions~\cite{wang_localgradientestimateharmonic_2010} when the Ricci curvature is quadratically asymptotically nonnegative (see also~\cite[Proposition 2.27]{benatti_asymptoticbehaviourcapacitarypotentials_2022}). However, as the proof presented in~\cite{wang_localgradientestimateharmonic_2010} is purely integral, integrating by parts the term containing the Ricci curvature and exploiting the $\CS^1$ decay of the metric leads to the following.

\begin{Lemma}\label{prop:gradientdecay}
 Let $(M, g)$ be a $\CS^1$-asymptotically flat Riemannian $3$-manifold with compact possibly empty boundary. Let $\Omega \subset M$ be a closed bounded subset with $\CS^1$-homologous to $\partial M$, and let $u_p$ be the solution to~\eqref{eq:p-cap-potential}, $1< p < 3$. Then, there exist $\kst >0$ and $R>0$ such that
 \begin{equation}
\label{eq:gradientdecay}
\abs{\D u_p}(x) \leq \kst\frac{u_p(x)}{\abs{x}}
 \end{equation}
 on $\set{\abs{x} \geq R}$.
\end{Lemma}
\begin{proof} We drop the subscript $p$. Let $\abs{x} \geq R$ for some $R$ large enough so that $\abs{\Gamma^k_{ij}} \leq {\kst}/\abs{x}$ $ \abs{g_{ij}-\delta_{ij}}\leq \kst$ on $\set{\abs{x} \geq R/3}$ for some $\kst>0$ and every $i,j,k=1,2,3$ and $B_{\abs{x}/2}(x) \subseteq\set{\abs{x} \geq R/3}$ for all $x \in \set{\abs{x} \geq R}$. The constant $\kst$ may change during the proof, but its value depends only on $R$, $g$ and $p$, not on $u$ or $x$. We explain how to modify the proof of the Cheng--Yau inequality in the ball $B=B_{\abs{x}/2}(x)$ presented in~\cite{wang_localgradientestimateharmonic_2010} to replace their lower bound on the Ricci curvature with the $\CS^1$-decay of the metric coefficients. The proof begins with an integral version of the Bochner formula
\begin{align}
\int_{B} \lop (f) \psi \dif\mu ={}& 2 \int_{B} f^{\frac{p}{2} - 1} \abs{\D\D u}^2 \psi \dif\mu + \bigg(\frac{p}{2} -1\bigg)\int_{B}\abs{\D f}^2 f^{\frac{p}{2} - 2} \psi \dif \mu \nonumber \\
&{}+\int_{B}f^{\frac{p}{2} - 1} \Ric(\D u, \D u) \psi \dif\mu \label{eq:primopassaggioyau}
\end{align}
for any $\psi \in \CS^\infty_c (B)$, where $f=\abs{\D u}^2$ and
\[
 \lop(f) = \div \bigg[ f^{p/2-1} \bigg( \D f + (p-2) \ip{\D u| \D f}\frac{ \D u }{\abs{ \D u }^2} \bigg) \bigg] - p f^{p/2-1} \ip{ \D u | \D f}.
\]
The only term containing the second derivatives of the metric is the one involving the Ricci curvature tensor. It can be written as
\begin{align}\label{eq:cheng-yau_Ricci_step1}
 \int_{B}\!f^{\frac{p}{2} - 1}\!\Ric(\D u, \D u) \psi \dif\mu = \! \int_{B} \!\! \big(\partial_k \Gamma^k_{ij} - \partial_i \Gamma^k_{kj} + \Gamma^k_{ij}\Gamma^m_{km} - \Gamma^m_{ik} \Gamma^k_{jm}\big)\D^i u\,\D^ju\, f^{\frac{p}{2} - 1}\psi \dif\mu. \!
\end{align}
The terms containing products of Christoffel symbols can be estimated using the $\CS^1$-asymptotic behaviour of the metric, using the fact we are on $\set{\abs{x}\geq R}$. Indeed,
\begin{align}\label{eq:cheng-yau_Ricci_step2}
 \int_{B}\big(\Gamma^k_{ij}\Gamma^m_{km} - \Gamma^m_{ik} \Gamma^k_{jm}\big)\D^i u\,\D^ju\,f^{\frac{p}{2} - 1}\psi \dif\mu\geq - \frac{\kst}{\abs{x}^2} \int_B f^{\frac{p}{2}} \psi \dif \mu.
\end{align}
On the other hand, using integration by parts, we have
\begin{align}
 \int_{B}\partial_k \Gamma^k_{ij}\D^i u\,\D^ju\, f^{\frac{p}{2} - 1}\psi \dif\mu &{}=- \int_{B}\Gamma^k_{ij}\partial_k \big(\D^i u\,\D^ju\, f^{\frac{p}{2} - 1}\psi \sqrt{ \det g}\big)\dif\mu_\delta \nonumber\\
 &{}\geq - \kst \int_B \frac{1}{\abs{x}}\abs{\D\D u}\abs{\D u} f^{\frac{p}{2}-1} \psi + f^{\frac{p}{2}}\bigg( \frac{\psi}{\abs{x}^2} + \frac{\abs{\D \psi}}{\abs{x}}\bigg)\!\dif \mu \nonumber\\
 &{} \geq - \kst \int_B \varepsilon^2\abs{\D\D u}^2f^{\frac{p}{2}-1} \psi +f^{\frac{p}{2}}\left( \frac{\big(1+\varepsilon^2\big)\psi}{\varepsilon^2\abs{x}^2} + \frac{\abs{\D \psi}}{\abs{x}}\right) \!\dif \mu \label{eq:cheng-yau_Ricci_step3}
\end{align}
for any $\varepsilon>0$, where we employed Young's inequality in the last step to the product between $\big(\abs{x}^{-1}\abs{\D\D u}\big)$ and $\abs{\D u}$, recalling that $f = \abs{\D u}^2$. Observe that we used the $\CS^1$-asymptotic behaviour of the metric not only to control $\abs{\partial g}$, but also to estimate $\abs{\partial \partial u}^2$ in terms of $\abs{\D \D u}^2$ and $\abs{x}^{-2}\abs{\D u}^2$. We can deal with the remaining term in the same way. Combining~\eqref{eq:cheng-yau_Ricci_step1}--\eqref{eq:cheng-yau_Ricci_step3}, we finally get
\begin{equation}\label{eq:cheng-yau_Ricci_final}
 \int_{B}f^{\frac{p}{2} - 1} \Ric(\D u, \D u) \psi \dif\mu \geq - \kst \int_B \varepsilon^2\abs{\D\D u}^2f^{\frac{p}{2}-1} \psi +f^{\frac{p}{2}}\left( \frac{\big(1+\varepsilon^2\big)\psi}{\varepsilon^2\abs{x}^2} + \frac{\abs{\D \psi}}{\abs{x}}\right) \dif \mu
\end{equation}
Following the argument in~\cite[Theorem 1.1]{wang_localgradientestimateharmonic_2010}, one can chose $\psi = f^b \eta^2$ for some $b > 1$ and with $\abs{\D \eta} \leq \mathrm{C} \, \eta /\abs{x}$. With this specification, the integrand of the last term in~\eqref{eq:cheng-yau_Ricci_final} is pointwise estimated by
\[
\frac{\abs{\D \psi}}{\abs{x}} f^{\frac{p}{2}} \leq \kst \psi\left(\varepsilon^2 \abs{\D \D u}^2 f^{\frac{p}{2} - 1}+ \frac{\big(1+\varepsilon^2\big)}{\varepsilon^2\abs{x}^2} f^{\frac{p}{2}} \psi\right).
\]
The last term in~\eqref{eq:cheng-yau_Ricci_final} can be absorbed in the others. Plugging this into~\eqref{eq:primopassaggioyau}, we deduce
\begin{align}
\int_{B} \mathscr{L} (f) \psi \dif\mu \geq{}& \big(2 - \epsilon^2\big)\int_{B}f^{\frac{p}{2} - 1} \abs{\D\D u}^2 \psi \dif\mu + \left(\frac{p}{2} -1\right)\int_{B}\abs{\D f}^2 f^{\frac{p}{2} - 2} \psi \dif \mu \nonumber\\
&{}- \frac{\kst\big(1+\varepsilon^2\big)}{\varepsilon^2\abs{x}^2} \int_{B}f^{\frac{p}{2}}\psi \dif\mu \label{eq:quasiultimopassaggioyau}
\end{align}
for any $\epsilon > 0$. Plugging now the last displayed identity at the bottom of~\cite[p. 763]{wang_localgradientestimateharmonic_2010} into~\eqref{eq:quasiultimopassaggioyau}, we get, choosing $\epsilon > 0$ small enough (depending only on~$p$), the inequality \cite[formula~(2.3)]{wang_localgradientestimateharmonic_2010}, with $\kappa$ given by a suitable uniform constant multiplying $\abs{x}^{-2}$. From this point on, the proof can be followed line by line, and yields the Cheng--Yau inequality of~\cite[Theorem~1.1]{wang_localgradientestimateharmonic_2010} in terms of $\kappa$ above in the ball $B_{\abs{x}/4}(x)$. This is exactly the claimed~\eqref{eq:gradientdecay}. \end{proof}

\begin{Remark}\label{rmk:chengwp}
We can rewrite the above estimate in terms of $w_p = -(p-1)\log u_p$. It reads
\[
\abs{\D w_p(x)} \leq \frac{ \kst}{\abs{x}}
\]
on $\set{\abs{x} \geq R}$, for $R$ large enough and some positive constant $\kst>0$ depending on $p$.
\end{Remark}
The following is a double-sided control on the solution $u_p$ to~\eqref{eq:p-cap-potential} of a bounded $\Omega \subset M$ with smooth boundary with respect to the Euclidean distance.

\begin{Lemma}\label{prop:menga_nomore}
 Let $(M, g)$ be a $\CS^1$-asymptotically flat Riemannian $3$-manifold with compact possibly empty boundary. Let $\Omega \subset M$ be a closed bounded subset with $\CS^1$-homologous to $\partial M$, and let $u_p$ be the solution to~\eqref{eq:p-cap-potential}, $1< p < 3$. Then, there exist $\kst >0$ and $R>0$ such that
 \[
 \kst^{-1}\abs{x}^{-\frac{3-p}{p-1}} \leq u_p(x) \leq \kst\abs{x}^{-\frac{3-p}{p-1}}
 \]
on $\set{\abs{x} \geq R}$.
\end{Lemma}

\begin{proof} We drop subscript $p$. The rightmost inequality follows by~\cite[Theorem 3.6]{mari_flowlaplaceapproximationnew_2022}, since having a positive isoperimetric constant is equivalent to having a global Sobolev inequality.

We get the leftmost inequality adapting the argument used for~\cite[Proposition 5.9]{holopainen_volumegrowthgreenfunctions_1999}. Integrating Lemma~\ref{prop:gradientdecay} as in~\cite{wang_localgradientestimateharmonic_2010} we have a Harnack inequality holding on large coordinate spheres
\[
\max_{\set{\abs{x} =r}} u \leq \kst \min_{\set{\abs{x} =r}} u_p,
\]
where $\kst$ does not depend on $r$. We are now committed to proving that
\[
 \max_{\set{\abs{x} =r}} u \geq \kst r^{-\frac{3-p}{p-1}},
\]
which concludes the proof. Let $m = \max\set{ u(x) \st \abs{x}=r}$. Then
\[
 \ncapa_p(\set{\abs{x} \leq r}) \geq \ncapa_p(\set{ u \geq m})= m^{-(p-1)} \ncapa_p(\partial \Omega).
\]
Using~\cite[Theorem 2.6]{heinonen_nonlinearpotentialtheorydegenerate_2018}, we have
\begin{align}
\label{eq:passaggio-pcapcondensatore}
 m\,\ncapa_p(\partial \Omega)^{-\frac{1}{p-1}} &\geq \ncapa_p(\set{\abs{x} \leq r})^{-\frac{1}{p-1}}\geq \sum_{j=0}^{+\infty} \big( \ncapa_p\big(\set{\abs{x} \leq 2^{j}r}, \set{ \abs{ x} < 2^{j+1} r}\big)\big)^{-\frac{1}{p-1}},
\end{align}
where we point out that the $p$-capacity of a condenser $(K,A)$ where $A$ is a open subset of $M$ and $K$ is a compact subset of $A$, is just defined as
\[
 \ncapa_p(K, A) = \inf \set{\frac{1}{4\pi}\bigg(\frac{p-1}{3-p}\bigg)^{p-1}\int_{A\smallsetminus K} \abs{\D v}^p \dif \mu\st v \in \CS_c^\infty(A),\, v \geq 1\text{ on } K}.
\]
Picking now a test function $\psi_r\in \CS_c^\infty\big(\set{\abs{x} \leq 2^{j+1}r}\big)$ taking the value $1$ on $\set{\abs{x} \leq 2^{j}r}$ and such that $\abs{\nabla \psi_r} \leq \kst/r$, we directly estimate the right-hand side of~\eqref{eq:passaggio-pcapcondensatore} with
\begin{gather}
 \sum_{j=0}^{+\infty} \big( \ncapa_p\big(\set{\abs{x} \leq 2^{j}r}, \set{ \abs{ x} < 2^{j+1} r}\big)\big)^{-\frac{1}{p-1}}\nonumber\\
 \qquad{}\geq\frac{p-1}{3-p}(4 \pi)^{-\frac{1}{p-1}} \kst \sum_{j=0}^{+\infty} \bigg( \int_{\set{\abs{x}< 2^{j+1}r}} \abs{\nabla \psi_r}^p\dif \mu\bigg)^{-\frac{1}{p-1}}\nonumber\\
 \qquad{}\geq\frac{p-1}{3-p}(4 \pi)^{-\frac{1}{p-1}} \kst \sum_{j=0}^{+\infty} \bigg( \frac{(2^j r)^p}{\abs{ \set{\abs{x}\leq 2^{j+1}r}}}\bigg)^{\frac{1}{p-1}}\nonumber\\
 \qquad{}\geq\frac{p-1}{3-p}(4 \pi)^{-\frac{1}{p-1}}\kst \int_{2r}^{+\infty} \bigg( \frac{t}{\abs{\set{\abs{x} \leq t}}}\bigg)^{\frac{1}{p-1}} \dif t.\label{eq:CONDENSATOR-LIKE-TERMINATOR}
\end{gather}
Finally, since $\abs{\set{\abs{x}\leq t}}= t^3(4\pi/3+o(1))$ as $t\to +\infty$, we can choose $R$ such that $\abs{\set{\abs{x}\leq r}}\leq\kst r^3$ for every $r \geq R$, so that, by~\eqref{eq:passaggio-pcapcondensatore} and~\eqref{eq:CONDENSATOR-LIKE-TERMINATOR}, we get
\[
 m \geq \kst\int_{2r}^{+\infty} \bigg( \frac{t}{\abs{\set{\abs{x} \leq t}}}\bigg)^{\frac{1}{p-1}} \dif t \geq \kst r^{-\frac{3-p}{p-1}},
\]
which concludes the proof. \end{proof}

\begin{Corollary}
 We can rewrite the above estimates in terms of $w_p = -(p-1)\log u_p$. They read
 \[
 (3-p) \log \abs{x} -\kst^{-1} \leq w_p \leq (3-p) \log \abs{x} +\kst
 \]
 on $\set{\abs{x} \geq R}$, for $R$ large enough and some positive constant $\kst>0$ depending on $p$.
\end{Corollary}

 We conclude by resuming some basic asymptotic expansions for $w_p$, substantially worked out in~\cite{benatti_asymptoticbehaviourcapacitarypotentials_2022}. We specialise in the case of $\CS^1$-asymptotically flat Riemannian $3$-manifolds and take advantage of the above observations in order to get rid of any Ricci curvature assumption.
\begin{Lemma}\label{lem:asylemma}
 Let $(M, g)$ be a $\CS^1$-asymptotically flat Riemannian $3$-manifold with compact, possibly empty boundary. Fix $1<p<3$. Let $\Omega\subseteq M$ be a closed bounded subset with $\CS^1$-boundary homologous to $\partial M$, $w_p$ the solution to~\eqref{eq:pb-intro} starting at $\Omega$ and $\Omega_t= \set{w_p \leq t}$. Then, for every $1<q<3$
 \begin{enumerate}\itemsep=0pt
 \item[$(1)$]
 $ w_p = (3-p) \log \abs{x} -\log(\ncapa_p(\partial \Omega))+o(1)$ as $\abs{x} \to +\infty$,

 \item[$(2)$]
 $ \D^iw_p = (3-p) \frac{x^i}{\abs{x}^2}(1+o(1))$ as $\abs{x} \to +\infty$,

 \item[$(3)$]
 $ \lim_{t \to +\infty} \ee^{-\frac{3-q}{3-p}t}\ncapa_q(\partial \Omega_t)=1$,

 \item[$(4)$] 
 $ \lim_{t \to +\infty} \ee^{-\frac{2}{3-p}t}\abs{ \partial \Omega_t}=4 \pi$.
 \end{enumerate}
\end{Lemma}
\begin{proof}Items (1) and (2) follow with the same strategy of~\cite[Theorem~3.1]{benatti_asymptoticbehaviourcapacitarypotentials_2022}, replacing~\cite[Corollary~2.25 and Proposition 2.27]{benatti_asymptoticbehaviourcapacitarypotentials_2022} with Lemmas~\ref{prop:menga_nomore} and~\ref{prop:gradientdecay}, respectively. Having (1), item~(3) follows at once. Indeed, for every $\varepsilon>0$
\[
\set{(3-p)\log\abs{x}\leq t-\varepsilon} \subset \Omega_t \subset \set{(3-p) \log \abs{x} \leq t+\varepsilon},
\]
for sufficiently large $t$. Hence, by monotonicity of the $q$-capacity we obtain
\[
\ncapa_q(\set{(3-p)\log\abs{x}= t-\varepsilon}) \leq\ncapa_q( \partial \Omega_t)\leq \ncapa_q(\set{(3-p) \log \abs{x} = t+\varepsilon}).
\]
Dividing both sides by $\ee^{-(3-q)t/(3-p)}$ and passing to the limit as $t \to +\infty$, in virtue of~\cite[Lemma~2.21]{benatti_asymptoticbehaviourcapacitarypotentials_2022} we get
\[
 \ee^{- \frac{ 3-q}{3-p}\varepsilon} \leq \lim_{t \to +\infty} \ee^{-\frac{3-q}{3-p}} \ncapa_q(\partial \Omega_t)\leq \ee^{\frac{3-q}{3-p}\varepsilon },
\]
from which we infer item~(3) sending $\varepsilon\to 0^+$. Item~(4) follows as~\cite[Proposition~3.4]{benatti_asymptoticbehaviourcapacitarypotentials_2022} replacing~\cite[Theorems 1.1 and~3.1]{benatti_asymptoticbehaviourcapacitarypotentials_2022} with items~(1) and (2), respectively. \end{proof}

\begin{Remark}
Up to this point, the content of this Section can be extended in any dimension~${n \geq 3}$ with obvious modifications. From now on we focus on dimension $3$ since the monotonicity formulas introduced are peculiar to this dimension.
\end{Remark}

\subsection{Concepts of mass in nonlinear potential theory}
The classical Hawking mass $\ma_H$, that is
\begin{equation}\label{eq:hawking-mass}
 \ma_{H}(\partial \Omega)= \frac{ \abs{ \partial \Omega}^{\frac{1}{2}}}{16 \pi^{\frac{3}{2}}} \bigg( 4 \pi - \int_{\partial \Omega}\frac{\H^2}{4} \dif \sigma \bigg),
\end{equation}
for $\Omega\subset M$, monotonically increases along the level sets of the weak IMCF~\cite{huisken_inversemeancurvatureflow_2001}. Such a property is clearly not preserved in general when one replaces the weak IMCF with solutions $w_p$ to~\eqref{eq:pb-intro}. For this reason, we consider a different family of quasi-local masses. We will call \emph{$p$-Hawking mass} the quantity
\begin{equation}\label{eq:p-Hawkingmass}
 \ma_H^{\tp}(\partial \Omega) = \frac{\ncapa_p(\partial \Omega)^{\frac{1}{3-p}}}{8 \pi}\Bigg[4\pi + \int_{\partial \Omega}\frac{\abs{\D w_p}^2}{(3-p)^2}\dif \sigma - \int_{\partial \Omega}\frac{\abs{ \D w_p} }{(3-p)}\H \dif \sigma \Bigg]
\end{equation}
for $\partial \Omega \in \CS^{1}$ with weak second fundamental form in $L^2(\partial \Omega)$. This should be thought of as a $p$-version of the classical Hawking mass. In fact, it is immediately seen that the $p$-Hawking mass formally converges to the Hawking mass as $p \to 1^+$, having in mind that along the weak IMCF $w_1$ we have $\abs{\D w_1} = \H$ and that the $p$-capacity of an outward minimizing set recovers the perimeter in such limit~\cite[Theorem 1.2]{fogagnolo_minimisinghullscapacityisoperimetric_2022}.
Crucially, as the Hawking mass is monotone along the weak IMCF~\cite[Geroch Monotonicity Formula 5.8]{huisken_inversemeancurvatureflow_2001},
 so the function $t \mapsto \ma_H^{\tp}(\partial \Omega_t)$, for $\Omega_t= \set{w_p \leq t}$, is monotone nondecreasing, as proven in~\cite{agostiniani_riemannianpenroseinequalitynonlinear_2022} (see actually~Appendix~\ref{app:monotonicities} for the \emph{full} monotonicity result).
\begin{Theorem}\label{thm:AMMO_monotonicity}
 Let $(M,g)$ be a strongly $p$-nonparabolic Riemannian $3$-manifold, ${1<p<3}$, with nonnegative scalar curvature and with connected, compact, possibly empty boundary. Assume that $H_2(M, \partial M;\Z) = \set{0}$. Let $\Omega\subseteq M$ be a closed bounded subset with connected $\CS^1$-boundary homologous to $\partial M$ and with second fundamental form $\h \in L^2(\partial \Omega)$, $w_p$ the solution to~\eqref{eq:pb-intro} starting at $\Omega$ and $\Omega_t= \set{w_p \leq t}$. The function $t\mapsto \ma_H^{\tp}(\partial \Omega_t)$ defined in~\eqref{eq:p-Hawkingmass} admits a~monotone nondecreasing $\BV_{\loc}(0,+\infty)$ representative and
 \begin{gather}
 \frac{\dd}{\dd t} \ma^{\tp}_H (\partial \Omega_t)=\frac{\ncapa_p(\partial \Omega_t)^{\frac{1}{3-p}}}{(3-p) 8 \pi} \Bigg( 4\pi -\int_{\partial \Omega_t}\frac{\mathrm{R}^\top}{2} \dif \sigma +\int_{\partial \Omega_t} \frac{\vert \mathring{\h} \vert^2}{2} + \frac{ \mathrm{R}}{2} + \frac{ \abs{ \D^\top \abs{ \D w_p} }^2}{\abs{ \D w_p}^2}\dif \sigma \nonumber \\
 \hphantom{\frac{\dd}{\dd t} \ma^{\tp}_H (\partial \Omega_t)=\frac{\ncapa_p(\partial \Omega_t)^{\frac{1}{3-p}}}{(3-p) 8 \pi} \bigg(}{}
 + \int_{\partial \Omega_t}\frac{5-p}{p-1}\bigg(\frac{ \abs{\D w_p}}{3-p} - \frac{ \H}{2}\bigg)^2\dif \sigma\Bigg)\label{eq:AMMO_monotonicity}
 \end{gather}
 holds at every $t$ regular for $w_p$.
\end{Theorem}
The $p$-Hawking mass has the very useful feature of dominating the Hawking mass times a~constant involving the global Sobolev constant of the underlying Riemannian manifold.
\begin{Lemma}\label{thm:hawking_p-hawking}
 Let $(M,g)$ be a strongly $p$-nonparabolic Riemannian $3$-manifold, $1<p<3$, with compact possibly empty boundary. Assume that the Sobolev constant $\kst_S$ of $(M, g)$ is positive. Then for every outward minimising $\Omega\subset M$ with $\CS^{1}$-boundary homologous to $\partial M$ with second fundamental form $\h \in L^2(\partial \Omega)$, we have
 \[
 \ma^{\tp}_{H}(\partial \Omega) \geq \Bigg( \frac{(3-p) (p-1)^{p-1}}{2^{2p-1}\pi^{\frac{p-1}{2}} p^{p}\kst_{S}^{\frac{3}{2}(p-1)}} \Bigg)^{\frac{1}{3-p}} \ma_{H}(\partial \Omega).
 \]
 \end{Lemma}

\begin{proof} Observe that
\[
 \int_{\partial \Omega}\frac{\abs{ \D w_p} }{(3-p)}\H \dif \sigma - \int_{\partial \Omega}\frac{\abs{\D w_p}^2}{(3-p)^2}\dif \sigma =\int_{\partial \Omega} \frac{\H^2}{4} \dif \sigma - \int_{\partial \Omega} \bigg( \frac{\H}{2} - \frac{ \abs{ \D w_p}}{3-p}\bigg)^2 \dif \sigma\leq \int_{\partial \Omega} \frac{ \H^2}{4} \dif \sigma.
\]
It is then enough to proceed as in the proof of~\cite[Theorem 1.3]{fogagnolo_minimisinghullscapacityisoperimetric_2022} to prove that
\[
 \ncapa_p (\partial \Omega)^{\frac{1}{3-p}} \geq \Bigg( \frac{(3-p) (p-1)^{p-1}}{4 \pi p^{p}\kst_{ S}^{\frac{3}{2}(p-1)}} \Bigg)^{\frac{1}{3-p}}\abs{ \partial \Omega}^{\frac{1}{2}}.
\tag*{\qed}
\]
\renewcommand{\qed}{}
\end{proof}

In~\cite[Theorem 2]{bray_capacitysurfacesmanifoldsnonnegative_2008} (see also the new proof proposed in~\cite{oronzio_admmassareacapacity_2022}) then extended to all $1<p<3$~\cite{xiao_harmoniccapacityasymptoticallyflat_2016}, the authors prove an upper bound for the capacity in terms of the area and the Willmore deficit. Observe that in their proof the asymptotically flat condition is assumed only to grant the existence of an IMCF starting at some $\Omega$. Here we assume the existence of the IMCF by requiring that $(M,g)$ is strongly $1$-nonparabolic (see~\cite{benatti_isoperimetricriemannianpenroseinequality_2022}).
We report here the statement, referring the reader to~\cite{bray_capacitysurfacesmanifoldsnonnegative_2008,xiao_harmoniccapacityasymptoticallyflat_2016} for the proof.

\begin{Proposition}[Nonlinear version of Bray--Miao's estimate]\label{thm:p-braymiao}
Let $(M,g)$ be a strongly $1$\-/nonparabolic Riemannian $3$-manifold with nonnegative scalar curvature and with connected, compact, possibly empty boundary. Assume that $H_2(M, \partial M; \Z) = \set{0}$. Let $\Omega\subset M$ be closed, bounded with connected $\CS^{1}$-boundary homologous to $\partial M$ with $\h\in L^2(\partial \Omega)$. Then
 \[
 \ncapa_p(\partial \Omega) \leq \bigg( \frac{\abs{\partial \Omega}}{4\pi}\bigg)^{\frac{3-p}{2}} \leftidx{_2}{F}{_1}\bigg(\frac{ 1}{2}, \frac{3-p}{p-1}, \frac{2}{p-1}; 1-\frac{1}{16\pi}\int_{\partial \Omega}\H^2 \dif \sigma\bigg)^{{-(p-1)}},
 \]
where $\leftidx{_2}{F}{_1}$ is the hypergeometric function.
\end{Proposition}
We recall that the hypergeometric function satisfies the following useful relation
\begin{equation}\label{eq:hypergeometric_relation}
 \leftidx{_2}{F}{_1}\bigg(\frac{ 1}{2}, \frac{3-p}{p-1}, \frac{2}{p-1}; \frac{2\ma}{t}\bigg)=\frac{3-p}{p-1}t^{\frac{3-p}{p-1}}\int_t^{+\infty}s^{\frac{2}{p-1}}\bigg(1- \frac{2\ma}{s}\bigg)^{-\frac{1}{2}} \dif s,
\end{equation}
where $1<p<3$ and $\ma\in \R$. Here, the integrand on right-hand side is, up to a scaling factor, the radial derivative of the rotationally symmetric $p$-capacitary potential of the horizon of the Schwarzschild of mass $\ma$.

Combining the previous proposition, the minimality of $\partial M$ and the isoperimetric Riemannian Penrose inequality~\cite[Theorem~1.3]{benatti_isoperimetricriemannianpenroseinequality_2022}, we obtain a sharp Penrose-type inequality for the $p$-capacity of the boundary, for every $1<p<3$.

\begin{Theorem}\label{thm:sharp-p-penrose}
 Let $(M,g)$ be a strongly $1$-nonparabolic Riemannian $3$-manifold with nonnegative scalar curvature and with smooth, compact, connected, minimal outermost boundary. Then, for every $1<p<3$ it holds
 \begin{equation}\label{eq:isoperimetric_penrose_pcap}
 \ncapa_p(\partial M)^{\frac{1}{3-p}}\leq 2\Bigg(\sqrt{\pi}\frac{\Gamma\big(\frac{2}{p-1}\big)}{\Gamma\big(\frac{2}{p-1}-\frac{1}{2}\big)}\Bigg)^{-\frac{p-1}{3-p}} \ma_{\iso},
 \end{equation}
 where $\Gamma$ is the gamma function. Moreover, the equality holds in~\eqref{eq:isoperimetric_penrose_pcap} if and only if $(M,g)$ is isometric to
 \[
 \Bigg(\R^n\smallsetminus \set{ \abs{x} < 2 \ma_{\iso}}, \bigg( 1+ \frac{ \ma_{\iso}}{2 \abs{x}}\bigg)^4\big(\eum_{ij}\dd x^i \otimes \dd x^j\big)\Bigg).
 \]
\end{Theorem}

\begin{proof} By~\cite[Lemma 2.8]{benatti_isoperimetricriemannianpenroseinequality_2022}, under the above assumptions we have $H_2(M, \partial M; \Z) = \set{0}$. By Proposition~\ref{thm:p-braymiao}, we have
\[
 \ncapa_p(\partial M )^{\frac{1}{3-p}}\leq 2\Bigg(\sqrt{\pi}\frac{\Gamma\big(\frac{2}{p-1}\big)}{\Gamma\big(\frac{2}{p-1}-\frac{1}{2}\big)}\Bigg)^{-\frac{p-1}{3-p}} \sqrt{\frac{\abs{\partial M}}{16 \pi}}.
\]
Then,~\eqref{eq:isoperimetric_penrose_pcap} follows from~\cite[Theorem 1.3]{benatti_isoperimetricriemannianpenroseinequality_2022}. The equality in~\eqref{eq:isoperimetric_penrose_pcap} implies the equality in~\cite[Theorem 1.3]{benatti_isoperimetricriemannianpenroseinequality_2022} yielding the rigidity statement. \end{proof}

This result has been provided in~\cite[Theorem 4]{bray_capacitysurfacesmanifoldsnonnegative_2008} for $p=2$ and~\cite[Theorem 1.1]{xiao_harmoniccapacityasymptoticallyflat_2016} for every $1<p<3$, in terms of the $\ADM$ mass and in the asymptotic flat regime. One can recover such formulation applying Theorem~\ref{thm:sharp-p-penrose} in conjunction with Theorem~\ref{thm:equivalence-masses}, proved below.

\section[p-isocapacitary Riemannian Penrose inequality]{$\boldsymbol{p}$-isocapacitary Riemannian Penrose inequality}
\label{sec:3}
In establishing the asymptotic comparison between the $p$-Hawking mass~\eqref{eq:p-Hawkingmass} and the $p$-iso\-capacitary mass, the quantity
\begin{equation}\label{eq:p-geroch_modified} \tilde{\ma}^{\tp}_H(\partial \Omega)= \frac{\ncapa_p(\partial \Omega)^{\frac{1}{3-p}}}{4\pi (3-p)}\Bigg( 4 \pi - \int_{\partial \Omega} \frac{\abs{ \D w_p}^2}{(3-p)^2} \dif \sigma\Bigg)
\end{equation}
will naturally appear. This quantity is closely related to the ones studied in~\cite{chan_monotonicitygreenfunctions_2022,munteanu_comparisontheorems3dmanifolds_2023} (see Theorem~\ref{thm:CCLT_monotonicity} and Remark~\ref{rmk:cclt_mon} for details). In the following lemma, we discuss the monotonicity properties of~\eqref{eq:p-geroch_modified} along the level sets of the solution $w_p$ to~\eqref{eq:pb-intro} and its relations with the $p$-Hawking mass~\eqref{eq:p-Hawkingmass}.

\begin{Lemma}\label{thm:inequality_modifiedgeroch_AMMO}
Let $(M,g)$ be a strongly $p$-nonparabolic Riemannian $3$-manifold, $1<p<3$, with nonnegative scalar curvature and with connected, compact, possibly empty boundary. Assume also that $(M,g)$ satisfies~\eqref{eq:energy_hp} and $H_2( M, \partial M; \Z)=\set{0}$. Let $\Omega\subseteq M$ be a closed bounded subset with connected $\CS^1$-boundary homologous to $\partial M$ and with second fundamental form $\h \in L^2(\partial \Omega)$, $w_p$ the solution to~\eqref{eq:pb-intro} starting at $\Omega$ and $\Omega_t= \set{w_p \leq t}$. Then the function $t \mapsto \tilde{\ma}^{\tp}_H(\partial \Omega_t)$ belongs to $W^{1,1}_{\loc}(0,+\infty)$ is monotone nondecreasing. Moreover, we have
\begin{equation}\label{eq:inequality_modifiedgeroch_AMMO}
 \ma^{\tp}_{H} ( \partial \Omega_t)\leq\tilde{\ma}^{\tp}_H(\partial \Omega_t)
\end{equation}
for every $t \in [0,+\infty)$, and
\begin{equation}\label{eq:limit_behaviour_phawking}
 \lim_{t \to+\infty} \ma^{\tp}_{H} ( \partial \Omega_t) = \lim_{t \to +\infty} \tilde{\ma}^{\tp}_H(\partial \Omega_t).
\end{equation}
\end{Lemma}

\begin{proof}
Denote
\begin{gather*}
 N(t) =\ncapa_p(\partial \Omega_t)^{-\frac{1}{p-1}}\Bigg( 4 \pi - \int_{\partial \Omega_t} \frac{\abs{ \D w_p}^2}{(3-p)^2} \dif \sigma\Bigg), \\
 D(t)= \ncapa_p(\partial \Omega_t)^{-\frac{2}{(3-p)(p-1)}},
\end{gather*}
we have that
\[
 \tilde{\ma}^{\tp}_H(\partial \Omega_t)=(\ncapa_p(\partial \Omega_t)^{\frac{1}{3-p}}\Bigg( 4 \pi - \int_{\partial \Omega_t} \frac{\abs{ \D w_p}^2}{(3-p)^2} \dif \sigma\Bigg) =\frac{N(t)}{D(t)}.
\]
Observe that, $N(t)$ is the quantity studied in~\cite{chan_monotonicitygreenfunctions_2022,munteanu_comparisontheorems3dmanifolds_2023}, while $1/D(t)$ is an exponentially growing term we multiplied it by. The function $N(t)$, and thus $\tilde{\ma}^{\tp}_H(\partial \Omega_t)=N(t)/D(t)$, belongs to~$W^{1,1}_{\loc}(0,+\infty)$ by Theorem~\ref{thm:CCLT_monotonicity}. Moreover,
\eqref{eq:weak_derivative} in Theorem~\ref{thm:CCLT_monotonicity} gives that $N'(t)/D'(t) = 4\pi (3-p)\ma_{H}^{\tp}(\partial \Omega_t)$ which is nondecreasing by Theorem~\ref{thm:AMMO_monotonicity}. Finally, $D(t)\to 0$ as $t \to +\infty$, while $N(t) \to 0$ as $t \to +\infty$ by the assumption \eqref{eq:energy_hp}. It is now a general fact that given two functions $f,g\in W^{1,1}_{\loc}(0,+\infty)$, with $g(t)\neq 0$ and $g'(t)\neq 0$, such that $f(t),g(t)\to 0$ as $t \to +\infty$ and with $f'(t)/g'(t)$ monotone nondecreasing, the function $f(t)/g(t)$ is monotone nondecreasing as well (see, e.g.,~\cite[Lemma 3.2]{zhu_comparisongeometryriccicurvature_1997}). Applying it with $f(t)=N(t)$ and $g(t)=D(t)$, we have that the function in $t \mapsto \tilde{\ma}^{\tp}_H(\partial \Omega_t)$ is nondecreasing. To prove~\eqref{eq:inequality_modifiedgeroch_AMMO}, it is enough to observe that
\[
 0\leq \frac{ \dd }{\dd t} \left( \frac{N(t)}{D(t)}\right) = \frac{N'(t)D(t) - N(t)D'(t)}{D(t)^2} = \frac{8 \pi}{(p-1)}\bigl(-\ma^{\tp}_H(\partial \Omega_t) + \tilde{\ma}^{\tp}_{H}(\partial \Omega_t)\bigr)
\]
for almost any $t$.
It then remains to prove~\eqref{eq:limit_behaviour_phawking}. On the other hand, by de L'H\^opital rule (see~\cite[Theorem A.1]{benatti_isoperimetricriemannianpenroseinequality_2022})
\[
 \lim_{t \to +\infty} 4\pi (3-p) \tilde{\ma}^{\tp}_H(\partial \Omega_t)= \lim_{t \to +\infty} \frac{N(t)}{D(t)} \leq \lim_{t \to +\infty} \frac{N'(t)}{D'(t)} = \lim_{t \to +\infty} 4 \pi (3-p) \ma_H^{\tp}(\partial \Omega_t).
\]
The reverse inequality easily follows by~\eqref{eq:inequality_modifiedgeroch_AMMO}. \end{proof}

The following result gives the $p$-capacitary counterpart of~\cite[Lemma 2.7]{benatti_isoperimetricriemannianpenroseinequality_2022}, that asymptotically controls the Hawking mass with the quasi-local isoperimetric mass of the evolving sets.

\begin{Lemma}[asymptotic comparison lemma]\label{thm:mass_control_pcap}
Let $(M,g)$ be a strongly $p$-nonparabolic Riemannian $3$-manifold, $1<p<3$, with nonnegative scalar curvature and with connected, compact, possibly empty boundary. Assume that $H_2( M, \partial M; \Z)=\set{0}$ and $(M,g)$ satisfies~\eqref{eq:energy_hp}. Let $\Omega\subseteq M$ be a closed bounded subset with connected $\CS^1$-boundary homologous to $\partial M$ and with second fundamental form $\h \in L^2(\partial \Omega)$, $w_p$ the solution to~\eqref{eq:pb-intro} starting at $\Omega$ and $\Omega_t= \set{w_p \leq t}$. Then
\begin{equation}\label{eq:mass_contro_pcap}
 \lim_{t \to +\infty} \ma_{H}^{\tp}(\partial \Omega_t)= \lim_{t \to +\infty} \tilde{\ma}_{H}^{\tp}(\partial \Omega_t) \leq\liminf_{t \to +\infty} \ma_{\iso}^{\tp}(\Omega_t),
\end{equation}
where $\Omega_t= \set{ w_p \leq t}$.
\end{Lemma}

\begin{proof} Assume that the right-hand side of~\eqref{eq:mass_contro_pcap} is finite, otherwise there is nothing to prove. The function $t\mapsto \abs{\Omega_t\smallsetminus \Crit w_p}$ is monotone continuous in $[0, +\infty)$, hence it is absolutely continuous. The generalised de L'H\^opital rule gives
\begin{align}
 \liminf_{t \to +\infty} \ma_{\iso}^{\tp}(\Omega_t)&{}\geq\liminf_{t \to +\infty} \frac{1}{2p\pi \ncapa_p(\partial \Omega_t)^{\frac{2}{3-p}}}\bigg( \abs{\Omega_t\smallsetminus \Crit w_p } -\frac{4 \pi}{3} \ncapa_p (\partial \Omega_t)^\frac{3}{3-p}\bigg)\nonumber\\
 &{}\geq \liminf_{t \to +\infty} \frac{(3-p)}{4 p\pi \ncapa_p(\partial \Omega_t)^{\frac{2}{3-p}}}\bigg(\int_{\partial \Omega_t} \frac{1}{\abs{\D w_p}} \dif \sigma-\frac{4 \pi}{3-p} \ncapa_p (\partial \Omega_t)^\frac{3}{3-p}\bigg).\label{eq:potential-delhopital}
\end{align}
By H\"{o}lder inequality, we have that
\[
 \int_{\partial \Omega_t }\frac{1}{\abs{\D w_p}} \dif \sigma \geq \bigg(\int_{\partial \Omega_t}\abs{\D w_p}^2 \dif \sigma\bigg)^{-\frac{p}{3-p}}\big( 4 \pi (3-p)^{p-1}\ncapa_p(\partial \Omega_t)\big)^{\frac{3}{3-p}}.
\]
Plugging it in~\eqref{eq:potential-delhopital}, we obtain
\[
 \liminf_{t \to +\infty} \ma_{\iso}^{\tp}(\Omega_t)\geq\liminf_{t \to +\infty} \frac{\ncapa_p(\partial \Omega_t)^{\frac{1}{3-p}}}{p\Big(\int_{\partial \Omega_t} \frac{\abs{ \D w_p}^2}{(3-p)^2} \dif \sigma\Big)^{\frac{p}{3-p}} }\left[ ( 4\pi)^{\frac{p}{3-p}}-\Bigg(\int_{\partial \Omega_t}\frac{\abs{ \D w_p}^2}{(3-p)^2} \dif \sigma \Bigg)^{\frac{p}{3-p}}\right].
\]
To simplify the notation, denote $f(z)=z^{p/(3-p)}$ and $z(t)=\int_{\partial \Omega_t} \abs{\D w_p}^2 /(3-p)^2 \dif \sigma$. Since $z(t)\leq 4 \pi$ in virtue of our assumptions, by Lemma~\ref{thm:inequality_modifiedgeroch_AMMO}
\begin{align}
 \liminf_{t \to +\infty} \ma_{\iso}^{\tp}(\Omega_t)&{}\geq \liminf_{t \to +\infty}\frac{1}{p f(z(t))} \frac{f(4 \pi )-f(z(t))}{4\pi-z(t)} \ncapa_p(\partial \Omega_t)^{\frac{1}{3-p}}( 4 \pi -z(t)) \nonumber\\
 &{}=\liminf_{t \to +\infty}\frac{4\pi(3-p)}{p f(z(t))} \frac{f(4 \pi )-f(z(t))}{4\pi-z(t)} \tilde{\ma}_{H}^{\tp}(\partial \Omega_t) \nonumber\\
 &{}\geq \liminf_{t \to +\infty}\frac{4\pi(3-p)}{p f(z(t))} \frac{f(4 \pi )-f(z(t))}{4\pi-z(t)} \ma_{H}^{\tp}(\partial \Omega_t).\label{eq:function_NPT_nonnegative}
\end{align}
The theorem follows once we prove the following claim.
\begin{Claim}
There exists a divergent increasing sequence $(t_n)_{n \in \N}$ realising the rightmost limit inferior of~\eqref{eq:function_NPT_nonnegative} and such that $z(t_n)\to 4 \pi$ as $n\to +\infty$.
\end{Claim}

Indeed, we would have
\[
 \lim_{n \to +\infty} \frac{4\pi}{f(z(t_n))} =(4 \pi )^{\frac{3-2p}{3-p}},\qquad \lim_{n\to +\infty}\frac{f(4 \pi )-f(z(t_n))}{4\pi-z(t_n)} = f'(4 \pi ) = \frac{p}{3-p} (4\pi)^{\frac{2p-3}{3-p}},
\]
that plugged into~\eqref{eq:function_NPT_nonnegative}, gives~\eqref{eq:mass_contro_pcap} in virtue of Theorem~\ref{thm:AMMO_monotonicity} and Lemma~\ref{thm:inequality_modifiedgeroch_AMMO}.

Let $t_n$ be a divergent increasing sequence $(t_n)_{n \in \N}$ realising the rightmost limit inferior of~\eqref{eq:function_NPT_nonnegative}. By Lemma~\ref{thm:inequality_modifiedgeroch_AMMO}, we have two possible cases:
\begin{enumerate}\itemsep=0pt
 \item[(1)] there exists $T>0$ such that $\tilde{\ma}^{\tp}_{H}(\partial \Omega_{t_n})\geq 0$ for all $t_n \geq T$, or
 \item[(2)] $\tilde{\ma}^{\tp}_{H}(\partial \Omega_{t_n})< 0$ for all $n\in \N$.
\end{enumerate}

\emph{Case~$1$.} Since $\tilde{\ma}^{\tp}_H(\partial \Omega_{t_n}) \geq 0 $, $z(t_n)\leq 4 \pi $ for every $t_n \geq T$. By contradiction, suppose there exists $\varepsilon>0$ such that $z(t_n) \leq 4 \pi - \varepsilon$ for every $n$ sufficiently large. Then, by~\eqref{eq:function_NPT_nonnegative}, there exists $\kst(p,\varepsilon)>0$ such that
\[
 +\infty>\liminf_{t\to +\infty} \ma^{\tp}_{\iso}(\Omega_{t}) \geq \lim_{n \to +\infty} \kst(p,\varepsilon)\ncapa_{p}(\partial \Omega_{t_n})^{\frac{1}{3-p}},
\]
which is clearly a contradiction. Hence, up to a not relabeled subsequence, $z(t_n) \to 4\pi$ as $n \to +\infty$. This proves the claim in this case.

\emph{Case~$2$.} 
Since $\tilde{\ma}^{\tp}_{H}(\partial \Omega_{t_n}) <0 $, $z(t_n) \geq 4 \pi$ for every $n\in \N$. Suppose by contradiction $z(t_n) \geq 4 \pi + \varepsilon$ for some $\varepsilon>0$. Then, by Theorem~\ref{thm:AMMO_monotonicity}, there exists $\kst(p,\varepsilon)>0$ such that
\[
 \tilde{\ma}^{(p)}_{H}(\partial \Omega) \leq \lim_{t \to +\infty} \tilde{\ma}^{\tp}_{H}(\partial \Omega_t) \leq -\kst(p,\varepsilon)\lim_{n \to +\infty}\ncapa_{p}(\partial \Omega_{t_n})^{\frac{1}{3-p}}=-\infty,
\]
which is a contradiction since $\abs{\D w_p} \in \CS^0(\partial \Omega)$, proving the claim also in this case.\end{proof}

Differently from~\cite[Lemma 2.7]{benatti_isoperimetricriemannianpenroseinequality_2022}, here we assumed~\eqref{eq:energy_hp}. We already mentioned in the Introduction that this condition is very mild. In the following remark, we better specify our assertion.

\begin{Remark}\label{rmk:geometric_energy}
First of all, observe that
\begin{equation}\label{eq:cheng_boundedness}
 \int_{\partial \Omega_t}\frac{\abs{\D w_p}^2}{(3-p)^2} \dif \sigma \leq 4 \pi \ee^{t}\ncapa_p(\partial \Omega)\sup_{\partial \Omega_t} \frac{\abs{ \D w_p}^{3-p}}{(3-p)^{3-p}}.
\end{equation}
If $\Ric(x) \geq - 2\kappa^2$ for some $\kappa \in \R$ and every $x \in M$, by~\cite[Theorem 1.1]{wang_localgradientestimateharmonic_2010} we have that $\abs{ \D w_p} \leq \kst_1$ for some constant depending on $p$ and $\kappa$. In particular, for $1<p<2$,~\eqref{eq:energy_hp} is fulfilled. The case $p=2$ may be treated as in~\cite[Corollary 1.1]{chan_monotonicitygreenfunctions_2022}. For the same reason, if~$(M,g)$ is $\CS^1$-asymptotically flat Riemannian manifold~\eqref{eq:energy_hp} is implied for every $1<p\leq 2$ by Lemma~\ref{prop:gradientdecay} (see also Remark~\ref{rmk:chengwp}).

Alternatively, assuming that $(M,g)$ is $\CS^0$-asymptotically flat and the Ricci tensor satisfies $\Ric(x) \geq - 2 \kappa^2 /(1+ \dist(x,o))^2$ for some $\kappa \in \R$ a fixed $o\in M$ and for every $x\in M$, one can cover the whole range $1<p<3$. Indeed, by~\cite[Theorem 1.1]{wang_localgradientestimateharmonic_2010} and~\cite [Theorem 1.1]{benatti_asymptoticbehaviourcapacitarypotentials_2022}, $\abs{ \D w_p} \leq \kst_3 \ee^{-t}$ for some positive constat $\kst_3$ depending only on $p$, $\kappa$ and $\Omega$. Plugging it into~\eqref{eq:cheng_boundedness}, we infer that $\int_{\partial \Omega_t} \abs{ \D w_p}^2 \dif \sigma \leq \kst_4$ for a positive constant $\kst_4$.
\end{Remark}
We establish a nonsharp Penrose inequality for the $p$-isocapacitary mass in the generality of Theorem~\ref{thm:pmt-intro}.

\begin{proof}[Proof of Theorem~\ref{thm:pmt-intro}] Assume first that $\partial M = \varnothing$.
Then, we let $w_p = -(p-1)G_p$, where~$G_p$ is the $p$-Green function issuing from some point $o \in M$. Then, all of results stated above for~$w_p$ starting from a given set can be obtained with no modifications in the proofs for such limit case. Moreover, by the asymptotic development of the $p$-Green function at the pole (see~\cite[Theorem~2.4]{mari_flowlaplaceapproximationnew_2022}), we have
\[
 \lim_{t \to -\infty}\ma^{\tp}_H ( \partial \Omega_t) = 0.
\]
 Applying Lemma~\ref{thm:mass_control_pcap} and Theorem~\ref{thm:AMMO_monotonicity}, we deduce that $\ma^{\tp}_{\iso} \geq 0$, as claimed.

 We do now treat the case $\partial M \neq \varnothing$. Let $w_p$ the solution to~\eqref{eq:pb-intro} and define $\Omega_t = \set{w_p \leq t}$. Observe that we can write
 \[
 2\ncapa_p(\partial M)^{\frac{1}{(3-p)}} = 2 \ma_{H}^{\tp}(\partial M)+(3-p)\tilde{\ma}_{H}^{\tp}(\partial M).
 \]
 Then, combining the monotonicity of $\ma_{H}^{\tp}$, that of $\tilde{\ma}_{H}^{\tp}$ following from~\eqref{eq:weak_derivative}, and the asymptotic comparison Lemma~\ref{thm:mass_control_pcap}, we obtain
\begin{align}
 2\ncapa_p(\partial M)^{\frac{1}{(3-p)}} &{}\leq \lim_{t \to +\infty} 2 \ma_{H}^{\tp}(\partial \Omega_t)+(3-p)\tilde{\ma}_{H}^{\tp}(\partial \Omega_t)\nonumber\\
 &{} \leq \liminf_{t \to +\infty} (5-p) \ma_{\iso}^{\tp}(\Omega_t) \leq (5-p) \ma^{\tp}_{\iso}.\label{eq:monotonicitiescombined}
\end{align}

Finally, we just have to discuss the equality case in the positive mass theorem. By the just prove Penrose-type inequality, $\partial M$ must be empty. Let then again, as above $w_p = -(p-1) \log G_p$. We deduce from the argument that yielded the positivity of the mass that $\ma^{\tp}_H$ must actually be constant. In particular, the right-hand side of~\eqref{eq:AMMO_monotonicity} constantly vanishes along the flow. The isometry with flat $\R^n$ then follows through very classical computations, that can be performed following the lines of~\cite[proof of Main Theorem~2]{huisken_inversemeancurvatureflow_2001}.\end{proof}
\begin{Remark}
We could have used the monotonicity of $\ma_{H}^{\tp}$ alone in~\eqref{eq:monotonicitiescombined}, instead of combining with the monotonicity of $\tilde{\ma}_{H}^{\tp}$. However, this would have led to the worse constant $2$ in the right-hand side of~\eqref{eq:p-penroseintro-crociata}.
\end{Remark}

\section{Proof of Theorem~\ref{thm:penroseadmintro}}
\pushQED{\qed}
The proof of Theorem~\ref{thm:penroseadmintro} follows from an asymptotic equivalence of $p$-Hawking masses. As one can expect, the $p$-Hawking mass has a better behaviour along the level set flow of $w_p$, which is the solution to~\eqref{eq:pb-intro}. But interestingly, under the right assumption on the asymptotic flatness, it tends to coincide with (the superior limit of) the Hawking mass on large sets. Moreover, it is asymptotically controlled by the other $q$-Hawking mass for $1<q<3$. Here we employ both the monotonicity of the mass $\ma^{\tp}_H$ and the better asymptotic behaviour of $\tilde{\ma}^{\tp}_H$ defined in~\eqref{eq:p-geroch_modified}. Indeed, we will use the latter one to ensure that $\ee^{-t/(3-p)} \ma^{\tp}_H(\partial \set{w_p \leq t})=o(1)$ as $t \to +\infty$, which permits to trigger the computations in~\cite[formula~(2.12)]{agostiniani_riemannianpenroseinequalitynonlinear_2022}.

\begin{Proposition}\label{thm:asymptotic_behaviour_phawking}
 Let $(M,g)$ be a $\CS^1$-asymptotically flat Riemannian $3$-manifold with nonnegative scalar curvature and with connected, compact, possibly empty boundary. Assume also that $H_2(M, \partial M; \Z)= \set{0}$. Fix $1<p<3$. Let $\Omega\subseteq M$ be homologous to $\partial M$ with connected $\CS^1$-boundary and $\h \in L^2(\partial \Omega)$, $w_p$ the solution to~\eqref{eq:pb-intro} starting at $\Omega$ and $\Omega_t= \set{w_p \leq t}$. Then%
 \begin{equation}
 \label{eq:hawkingpq}
 \lim_{t \to +\infty }\ma^{\tp}_H(\partial \Omega_t)=\limsup_{t \to +\infty} \ma_H(\partial \Omega_t) \leq \limsup_{t \to +\infty} \ma_H^{\tp[q]}(\partial \Omega_t)
 \end{equation}
 for every $1<q<3$.
\end{Proposition}
\begin{proof} The inequality appearing in~\eqref{eq:hawkingpq} is obtained arguing as in Lemma~\ref{thm:hawking_p-hawking}. Indeed, we get%
\begin{equation} \label{eq:upper_bound_hawking}
\limsup_{t \to +\infty} \ma_H(\partial \Omega_t) \leq \limsup_{t \to +\infty} \ncapa_q(\partial \Omega_t)^{-\frac{1}{3-q}}\sqrt{\frac{\abs{ \partial \Omega_t}}{4\pi}}\ma_H^{\tp[q]} (\partial \Omega_t)=\limsup_{p \to +\infty} \ma^{\tp[q]}_H(\partial \Omega_t),
\end{equation}
where the last identity follows by Lemma~\ref{lem:asylemma}\,(3)\,(4). In order to show the identity appearing in~\eqref{eq:hawkingpq}, we are thus left to show the inequality
\begin{equation}
\label{eq:claim-ineq}
 \lim_{t \to +\infty }\ma^{\tp}_H(\partial \Omega_t) \leq \limsup_{t \to +\infty} \ma_H(\partial \Omega_t),
\end{equation}
the reverse one consisting in~\eqref{eq:upper_bound_hawking} with $p = q$. To do so, we claim that
\begin{equation}\label{eq:vanishing_condition}
 \Bigg[4\pi + \int_{\partial \Omega}\frac{\abs{\D w_p}^2}{(3-p)^2}\dif \sigma - \int_{\partial \Omega}\frac{\abs{ \D w_p} }{(3-p)}\H \dif \sigma \Bigg]=o(1)
\end{equation}
as $t \to +\infty$. Indeed, if this happens, we can follow the chain of inequalities in~\cite[(2.12)]{agostiniani_riemannianpenroseinequalitynonlinear_2022} (see also~\cite[Theorem 4.11]{benatti_isoperimetricriemannianpenroseinequality_2022}) and obtain
\begin{equation*}
 \lim_{t \to +\infty} \ma^{\tp}_H(\partial \Omega_t) \leq \limsup_{t \to +\infty} \ncapa_p(\partial \Omega_t)^{\frac{1}{3-p}}\sqrt{\frac{4\pi}{\abs{\partial \Omega_t}}} \ma_H(\partial \Omega_t)= \limsup_{t \to +\infty} \ma_H(\partial \Omega_t),
\end{equation*}
where again we applied Lemma~\ref{lem:asylemma}\,(3)\,(4), proving~\eqref{eq:claim-ineq}.

We then proceed to prove~\eqref{eq:vanishing_condition}.
If $\ma^{\tp}_H(\partial \Omega_t) <0$
 for every $t\in [0,+\infty)$, arguing as in Case~2 of the proof of Lemma~\ref{thm:mass_control_pcap}, we deduce that~\eqref{eq:vanishing_condition} must hold. Otherwise, we would contradict the monotonicity formulas in Theorem~\ref{thm:AMMO_monotonicity}. Conversely, appealing again to the monotonicity formulas in Theorem~\ref{thm:AMMO_monotonicity}, $t \mapsto \ma^{\tp}_H(\partial \Omega_t)$ must be definitely nonnegative. Observe that by
Lemma~\ref{lem:asylemma}\,(1)\,(2), we have
\[
4 \pi - \int_{\partial \Omega_t} \frac{ \abs{ \D w_p}^2}{(3-p)^2} \dif \sigma = o(1)
\]
as $t\to +\infty$. Hence, Lemma~\ref{thm:inequality_modifiedgeroch_AMMO} implies
\[
0\leq \ma^{\tp}_H(\partial \Omega_t) \leq \tilde{\ma}^{\tp}_H(\partial \Omega_t) = o(\ee^{t})
\]
as $t\to +\infty$. Dividing both sides by $\ncapa_p(\partial \Omega_t)^{1/(3-p)}$ we get~\eqref{eq:vanishing_condition}.
 \end{proof}

\begin{proof}[Conclusion of the proof of Theorem~\ref{thm:penroseadmintro}.]
 Differently from the case $p=1$, corresponding to the classical Hawking mass, here we assume connectedness of the boundary of the manifold. In fact, it is not clear to us how to adapt the argument employed in~\cite[Section 6]{huisken_inversemeancurvatureflow_2001}, where the authors took advantage of the horizons being minimal and outward minimizing in order to prescribe a jump that maintains the monotonicity of the Hawking mass. The difficulties when dealing with the $p$-Hawking mass arise in connection with the gradient of $w_p$ appearing in its expression. Assuming $\partial M$ to be connected, we can consider the solution $w_p$ to~\eqref{eq:pb-intro} starting at $\Omega= \partial M$ and $\Omega_t= \set{w_p\leq t}$.
The boundary of $M$ being outermost implies that $H_2(M, \partial M; \Z)= \set{0}$ (see~\cite[Lemma 4.1]{huisken_inversemeancurvatureflow_2001}, or the alternative argument in the proof of~\cite[Lemma 2.8]{benatti_isoperimetricriemannianpenroseinequality_2022}).
Applying Proposition~\ref{thm:asymptotic_behaviour_phawking} for $q=2$, we have
\[
 \frac{ \ncapa_{p}(\partial M)^{\frac{1}{3-p}}}{2} \leq \lim_{t \to +\infty} \ma^{\tp}_H(\partial \Omega_t) \leq \limsup_{t \to +\infty} \ma^{\tp[2]}_H(\partial \Omega_t).
\]
Since by Lemma~\ref{lem:asylemma}\,(2), $\partial \Omega_t$ is regular for any $t$ large enough, we can use~\cite[Theorem~4.11]{benatti_isoperimetricriemannianpenroseinequality_2022} to control the right-hand side with $\ma_{\ADM}$, concluding the proof of~\eqref{eq:p-penroseintro}. Observe now that an outermost minimal boundary is outward minimising. If this were not the case, the outward minimising hull~\cite{fogagnolo_minimisinghullscapacityisoperimetric_2022, huisken_inversemeancurvatureflow_2001} would be a closed minimal surface homologous to it and, by the Maximum Principle, disjoint from $\partial M$. Then, letting $p \to 1^+$ and appealing to~\cite[Theorem~1.2]{fogagnolo_minimisinghullscapacityisoperimetric_2022} recovers the sharp Penrose inequality~\eqref{eq:penroseintro}.
\end{proof}

\section[Relation between the isoperimetric mass and the p-isocapacitary mass]{Relation between the isoperimetric mass and the \\ $\boldsymbol{p}$-isocapacitary mass}
\label{sec:4}
We now employ the explicit control of the Hawking mass in terms of the $p$-Hawking mass Lemma~\ref{thm:hawking_p-hawking} to produce an upper bound on the isoperimetric mass in terms of $p$-isocapacitary mass. This bound is not sharp but sharpens as $p\to 1^+$.
\begin{Lemma}\label{cor:lowerbound_miso}
 Let $(M,g)$ be a $\CS^0$-asymptotically flat Riemannian $3$-manifold with nonnegative scalar curvature and with smooth, compact, minimal, possibly empty boundary. Assume that~$(M,g)$ satisfies~\eqref{eq:energy_hp} for some $1<p<3$. Then
 \begin{equation*}
 \ma_{\iso} \leq \Bigg( \frac{2^{2p-1}\pi^{\frac{p-1}{2}} p^{p}\kst_{ S}^{\frac{3}{2}(p-1)}}{(3-p) (p-1)^{p-1}} \Bigg)^{\mathrlap{\frac{1}{3-p}}}\;\;\; \ma^{\tp}_{\iso},
 \end{equation*}
 where $\kst_S$ is the global Sobolev constant of $(M,g)$.
\end{Lemma}

\begin{proof}
By the topological description of manifolds like these, reworked in~\cite[Lemma 2.8]{benatti_isoperimetricriemannianpenroseinequality_2022}, we can assume that our Riemannian manifold has a (possibly empty) minimal, outermost boundary such that $H_2(M,\partial M; \Z) = \set{0}$.
Let $E\subset M$ be a closed bounded subset with smooth boundary such that any connected component of $\partial M$ is either contained in $E$ or disjoint from $E$. Using~\cite[Theorem 6.1]{huisken_inversemeancurvatureflow_2001} (see also~\cite[Proposition 2.5]{benatti_isoperimetricriemannianpenroseinequality_2022}), we can find a subset $\Omega$ closed bounded with $\CS^{1}$-boundary homologous to $\partial M$ and with $\h \in L^2(\partial \Omega)$ such that $\ma_H(\partial E) \leq \ma_H(\partial \Omega)$. Let $w_p$ the solution to~\eqref{eq:pb-intro} starting at $\Omega$ and $\Omega_t = \set{ w_p \leq t}$ . By Lemmas~\ref{thm:hawking_p-hawking} and~\ref{thm:mass_control_pcap}, we now have
\[
\ma_H(\partial E) \leq \Bigg( \frac{2^{2p-1}\pi^{\frac{p-1}{2}} p^{p}\kst_{ S}^{\frac{3}{2}(p-1)}}{(3-p) (p-1)^{p-1}} \Bigg)^{\mathrlap{\frac{1}{3-p}}}\;\;\;\limsup_{t \to +\infty} \ma^{\tp}_H(\partial \Omega_t) \leq\Bigg( \frac{2^{2p-1}\pi^{\frac{p-1}{2}} p^{p}\kst_{ S}^{\frac{3}{2}(p-1)}}{(3-p) (p-1)^{p-1}} \Bigg)^{\mathrlap{\frac{1}{3-p}}}\;\;\;\ma^{\tp}_{\iso}.
\]
Since we have a control on the Hawking mass of every $E$, we can apply~\cite{jauregui_lowersemicontinuitymass0_2019} (see~\cite[Theorem~2.6]{benatti_isoperimetricriemannianpenroseinequality_2022} for the precise statement and remarks) to control the isoperimetric mass with the same quan\-tity.
\end{proof}

We prove a family of equivalent formulations for the $p$-isocapacitary masses, as well as for the isoperimetric one. The proofs will follow the one given in~\cite[Lemma 10]{jauregui_admmasscapacityvolumedeficit_2020} for the $2$-isocapacitary mass.
\begin{Proposition}\label{prop:equivalence_of_iso_p_mass}
Let $(M,g)$ be a $\CS^0$-asymptotically flat Riemannian $3$-manifold with compact, possibly empty boundary. Then, for $1<p<3$, we have
\begin{equation}\label{eq:equivalence_of_p_masses_for_different_exponents}
 \ma^{\tp}_{\iso}= \sup_{(\Omega_j)_{j \in \N}} \limsup_{j \to +\infty} \frac{2\ncapa_p(\partial \Omega)^{\frac{1-3 \alpha}{3-p}}}{3p\alpha} \bigg( \bigg(\frac{3\abs{ \Omega_j}}{4\pi}\bigg)^{\alpha} - \ncapa_p(\partial \Omega_j)^{\frac{3 \alpha}{3-p}}\bigg)
\end{equation}
for every $\alpha\geq 1/3$.
\end{Proposition}
The main computation performed in order to prove the result above is the following one, that we isolate for future reference.

\begin{Lemma}
 \label{lem:equivalence-specialexhaustion}
 Let $(M,g)$ be a $\CS^0$-asymptotically flat Riemannian $3$-manifold with compact possibly empty boundary, and $1<p<3$. Let $(\Omega_j)_{j \in \N}$ be an exhaustion of $M$ such that
 \begin{equation}
 \label{eq:limit_behaviour_on_goodsequence}
 \lim_{j \to +\infty} \frac{\abs{\Omega_j}}{\ncapa_p(\partial \Omega_j)^{\frac{3}{3-p}}} = \frac{4 \pi }{3}.
\end{equation}
Then, we have
\begin{equation*}
\limsup_{j \to + \infty}\ma^{\tp}_{\iso}(\Omega_j) = \limsup_{j \to +\infty} \frac{2\ncapa_p(\partial \Omega_j)^{\frac{1-3 \alpha}{3-p}}}{3p\alpha} \bigg( \bigg(\frac{3\abs{ \Omega_j}}{4\pi}\bigg)^{\alpha} - \ncapa_p(\partial \Omega_j)^{\frac{3 \alpha}{3-p}}\bigg)
\end{equation*}
for any $\alpha \in \R \setminus \set{0}$.
\end{Lemma}
\begin{proof}
 Let $(\Omega_j)_{j\in \N}$ be a sequence such that~\eqref{eq:limit_behaviour_on_goodsequence} holds. Up to considering a subsequence, we can assume that $(\Omega_j)_{j \in \N}$ realises the superior limit. Denote $f(z)= z^{\alpha}$, $z_j= \abs{\Omega_j}/\ncapa_p(\partial \Omega_j)^{3/(3-p)}$, we have that
\begin{equation}\label{eq:limit-p-iso}
 \limsup_{j \to +\infty} \ma^{\tp}_{\iso}(\Omega_j) = \lim_{j\to +\infty} \frac{\ncapa_p(\partial \Omega_j)^{\frac{1}{3-p}} }{2p \pi } \frac{ z_j - 4 \pi /3 }{f(z_j)- f( 4\pi/3)} (f(z_j) - f( 4 \pi /3)).
\end{equation}
Since $f$ is differentiable at $4\pi/3$ and $z_j\to 4 \pi/3\neq 0$ as $j\to +\infty$ by~\eqref{eq:limit_behaviour_on_goodsequence}, we have
\[
 \lim_{j \to + \infty} \frac{ z_j - 4 \pi /3 }{f(z_j)- f( 4\pi/3)} = \frac{1}{f'(4 \pi /3)}= \frac{3 ^{\alpha-1}}{\alpha(4 \pi )^{\alpha-1}}.
\]
Plugging this into~\eqref{eq:limit-p-iso} we conclude.
\end{proof}

\begin{proof}[Proof of Proposition~\ref{prop:equivalence_of_iso_p_mass}] We claim that it is enough to prove the equivalence on sequences such that
\eqref{eq:limit_behaviour_on_goodsequence} holds, so that Proposition~\ref{prop:equivalence_of_iso_p_mass} follows from Lemma~\ref{lem:equivalence-specialexhaustion}.
 Let then $(\Omega_j)_{j \in \N}$ be an exhaustion. By the $p$-isocapacitary inequality, we have that
\[
 \limsup_{j \to +\infty} \frac{\abs{\Omega_j}}{\ncapa_p(\partial \Omega_j)^{\frac{3}{3-p}}} \leq \frac{4 \pi }{3}.
\]
Indeed, the metric $g$ becomes uniformly equivalent to the flat Euclidean metric on $M \smallsetminus \Omega_j$ as $j \to +\infty$.
Moreover, for sufficiently large $j$ there exists a unique $r_j>0$ such that the coordinate ball $B_{r_j}$ has the same volume of $\Omega_j$. Define
\[
 \Omega'_j= \begin{cases}
 \Omega_j & \text{if } \capa_p(\Omega_j)\leq \capa_p(B_{r_j}),\\
 B_{r_j} & \text{if } \capa_p(\Omega_j)> \capa_p(B_{r_j}).
 \end{cases}
\]
The sequence $\big(\Omega'_j\big)_{j\in \N}$ is an exhaustion of $M$ and
\[
 \liminf_{j \to +\infty} \frac{\big|\Omega'_j\big|}{\ncapa_p(\partial \Omega'_j)^{\frac{3}{3-p}}}\geq\liminf_{j \to +\infty} \frac{\big|B_{r_j}\big|}{\ncapa_p\big(\partial B_{r_j}\big)^{\frac{3}{3-p}}}= \frac{4 \pi }{3},
\]
where the right-hand side is computed using the asymptotic flatness. In particular, the sequence~$(\Omega'_j)_{j\in \N}$ fulfils~\eqref{eq:limit_behaviour_on_goodsequence}, $\abs{\Omega'_j}= \abs{\Omega_j}$ and $\ncapa_p(\partial\Omega'_j)\leq \ncapa_p(\partial\Omega_j)$. Then, when $\alpha \geq 1/3$, $(\Omega'_j)_{j \in \N}$ is a better competitor both for $\ma^p_{\iso}$ as in the definition of $p$-isocapacitary mass~\eqref{eq:isopcapacitary_mass} and for the right-hand side of~\eqref{eq:equivalence_of_p_masses_for_different_exponents}. This completes the proof.
 \end{proof}

Completely analogous results hold for the perimeter and the isoperimetric mass. We gather them in the following statement.
\begin{Proposition}
 Let $(M,g)$ be a $\CS^0$-asymptotically flat Riemannian $3$-manifold with compact possibly empty boundary. Let $(\Omega_j)_j$ be an exhaustion of $M$ such that
 \begin{equation*}
 \lim_{j\to +\infty} \frac{ \abs{ \Omega_j}}{\abs{\partial \Omega_j}^{\frac{3}{2}}} = \frac{1}{6 \sqrt{\pi}}.
 \end{equation*}
 Then
 \begin{equation}
 \label{eq:equivalencesmassesformula}
\limsup_{j \to +\infty}\frac{2}{\abs{\partial \Omega_j}}\bigg( \abs{\Omega_j} -\frac{\abs{\partial \Omega_j}^{\frac{3}{2}}}{6 \sqrt{\pi}}\bigg) = \limsup_{j\to +\infty}\frac{\abs{ \partial \Omega_j}^{\frac{1-3 \alpha}{2}}}{3 \alpha \sqrt{\pi}} \big( (6 \sqrt{\pi}\abs{ \Omega_j})^{\alpha} - \abs{ \partial \Omega_j}^{ \frac{3\alpha}{2}} \big)
\end{equation}
holds for any $\alpha \in \R\smallsetminus\set{0}$.
As a consequence, we have
 \begin{equation*}
 \ma_{\iso}= \sup_{(\Omega_j)_{j \in \N}} \limsup_{j \to +\infty}\frac{\abs{ \partial \Omega_j}^{\frac{1-3 \alpha}{2}}}{3 \alpha \sqrt{\pi}} \big( (6 \sqrt{\pi}\abs{ \Omega_j})^{\alpha} - \abs{ \partial \Omega_j}^{ \frac{3\alpha}{2}} \big)
 \end{equation*}
 for every $\alpha\geq 1/3$.
\end{Proposition}
The inequality $\ma^{\tp}_{\iso} \leq \ma_{\iso}$ will substantially be a consequence of the following $p$-isocapacitary inequality for sets with volume going to infinity. Its isoperimetric version was pointed out in~\cite[Corollary C.3]{chodosh_isoperimetryscalarcurvaturemass_2021}.
\begin{Theorem}[sharp asymptotic $p$-isocapacitary inequality]\label{thm:sharp-isopcapacitary}
 Let $(M,g)$ be a $\CS^0$-asymptotically flat Riemannian $3$-manifold with compact possibly empty boundary $\partial M$. Then, for every $1<p<3$, we have that \begin{equation}\label{eq:sharp_iso_p_capacitary_inequality}
 \abs{ \Omega}^{\frac{3-p}{3}}\leq \bigg( \frac{4 \pi }{3} \bigg)^{\frac{3-p}{3}} \ncapa_p(\partial \Omega)+\frac{p(3-p)}{2} \ma \bigg( \frac{4 \pi }{3} \bigg)^{\frac{3-p}{3}}\ncapa_p(\partial \Omega)^{\frac{2-p}{3-p}}(1+o(1))
 \end{equation}
 as $\abs{ \Omega} \to +\infty$, where $\Omega$ closed and bounded with $\CS^{1, \alpha}$-boundary containing $\partial M$ and $\ma>-\infty$ is such that $\ma\geq \ma_{\iso}$.
\end{Theorem}

\begin{proof} Assume that $\ma<+\infty$, otherwise there is nothing to prove. We claim that for every~${\varepsilon>0}$ small enough, there exists $V_\varepsilon>\varepsilon^{-3}$ such that
\begin{equation}\label{eq:sharp_isoperimetric_inequality}
 \big(6 \sqrt{\pi} \abs{\Omega}\big)^{\frac{2p}{3}} \leq \abs{\partial \Omega}^{p}+ 2p\sqrt{\pi}(\ma+\varepsilon)\abs{ \partial \Omega}^{\frac{2p-1}{2}}
\end{equation}
for every $\Omega \subseteq M$ such that $\abs{\Omega} \geq V_\varepsilon$. Indeed, if this were not the case, we would find a~sequence~$(\Omega_j)_{j \in \N}$ with $\abs{\Omega_j} \to + \infty$ such that the right-hand side with $\alpha=2p/3\geq 1/3$, and thus the left-hand side, of~\eqref{eq:equivalencesmassesformula}, is strictly bigger than $ \ma \geq \ma_{\iso}$. Since, by the isoperimetric inequality, the perimeters of the $\Omega_j$'s diverge at infinity too, this would contradict~\cite[Proposition 37]{jauregui_lowersemicontinuitymass0_2019}, stating that one can relax the competitors in the definition of the isoperimetric mass in order to include any sequence of bounded sets containing $\partial M$ with diverging perimeters.

We can now assume that
\begin{equation}\label{eq:trivial assumption}
 \abs{\Omega}^{\frac{3-p}{2p}} \geq \bigg(\frac{4 \pi}{3}\bigg)^{\frac{3-p}{2p}}\ncapa_p(\partial \Omega)^{\frac{3}{2p}}(1-\varepsilon)^{\frac{3-p}{2p}},
\end{equation}
otherwise~\eqref{eq:sharp_iso_p_capacitary_inequality} is trivially satisfied. Let $w_p\colon M \smallsetminus \Omega \to \R$ be the solution to~\eqref{eq:pb-intro} starting at $\Omega$, $w_p = -(p-1) \log u_p$ and let $\Omega_t=\set{u_p\geq t}\cup \Omega $ and $V(t)= \abs{ \Omega_t}\geq V_\varepsilon$ for every $t \in (0,1)$. The H\"{o}lder's inequality with exponents $a=p$ and $b=p/(p-1)$ gives
\begin{equation}\label{eq:level_area_bound_CS}
 \abs{ \partial \Omega_t}^p \leq \bigg(\,\int_{ \partial \Omega_t} \abs{\D u_p}^{p-1} \dif \sigma \bigg)\bigg(\, \int_{\partial \Omega_t}\frac{1}{| \D u_p |}\dif \sigma\bigg)^{p-1} = 4\pi \ncapa_p \bigg( \frac{3-p}{p-1}\bigg)^{p-1}[-V'(t)]^{p-1}
\end{equation}
for almost every $t \in (0,1]$, where $\ncapa_p=\ncapa_p(\partial \Omega)$. Plugging it into~\eqref{eq:sharp_isoperimetric_inequality} and integrating on $(0, 1)$ we obtain
\begin{gather}
 \int_0^1\frac{\big[6 \sqrt{\pi}V(t)\big]^{\frac{2p}{3}}}{(-V'(t))^{p-1}}\dif t \nonumber\\
\qquad{} \leq 4\pi \ncapa_p \bigg( \frac{3-p}{p-1}\bigg)^{p-1}+ \Bigg[4\pi\ncapa_p \bigg( \frac{3-p}{p-1}\bigg)^{p-1}\Bigg]^{\frac{2p-1}{2p}} \int_0^1\frac{2p\sqrt{\pi}(\ma+\varepsilon)}{(-V'(t))^{\frac{p-1}{2p}}}\dif t. \label{eq:integration_sharpcap}
\end{gather}
Applying~\eqref{eq:sharp_isoperimetric_inequality} and the isoperimetric inequality, we get
\[
 (6 \sqrt{\pi}V(t))^{\frac{2p}{3}} \leq \abs{\partial \Omega_t}^{p}\bigg(1+ \frac{2p\sqrt{\pi}\abs{\ma+\varepsilon}}{\sqrt{\abs{\partial \Omega_t}}}\bigg)\leq\abs{\partial \Omega_t}^{p}\bigg(1+ \frac{\kst}{V(t)^{\frac{1}{3}}}\bigg),
\]
where $\kst$ depends only on $\ma$, $p$ and the isoperimetric constant. Plugging it into~\eqref{eq:level_area_bound_CS}, we have that
\begin{equation}\label{eq:lower_bound_on_derivative}
 [-V'(t)]^{p-1} \geq \bigg( \frac{p-1}{3-p}\bigg)^{p-1} \frac{3^{\frac{2p}{3}} V(t)^{\frac{2p}{3}}}{(4\pi)^{\frac{3-p}{p}}\ncapa_p}\bigg(1+ \frac{\kst}{V(t)^{\frac{1}{3}}}\bigg)^{-1}.
\end{equation}
 Hence, using \eqref{eq:lower_bound_on_derivative}, the assumption $V(t)\geq \abs{\Omega}\geq \varepsilon^{-3}$, a change of variable in the integral and~\eqref{eq:trivial assumption} yield
\begin{align}
 \int_0^1\bigl[-V'(t)\bigr]^{-\frac{p-1}{2p}}\dif t&{}=-\int_0^1\bigl[-V'(t)\bigr]^{-\frac{3p-1}{2p}}V'(t)\dif t \nonumber\\
 &{} \leq-\int_0^1\Bigg[ \bigg( \frac{3-p}{p-1}\bigg)^{p-1} \frac{(4\pi)^{\frac{3-p}{p}}\ncapa_p}{3^{\frac{2p}{3}} V(t)^{\frac{2p}{3}}}\bigg(1+ \frac{\kst}{V(t)^{\frac{1}{3}}}\bigg)\Bigg]^{\frac{3p-1}{2p(p-1)}} V'(t) \dif t \nonumber\\
 &{}\leq\Bigg[\bigg(\frac{3-p}{p-1}\bigg)^{p-1}\frac{(4 \pi )^{\frac{3-p}{3}}}{3^{\frac{2p}{3}}}\ncapa_p(1+\kst \varepsilon)\Bigg]^{\frac{3p-1}{2p(p-1)}}\int_{\abs{\Omega}}^{+\infty} V^{-\frac{3p-1}{3(p-1)}} \dif V \nonumber\\
 &{}\leq\frac{3-p}{2}\left[\bigg(\frac{3-p}{p-1}\bigg)^{\frac{(p-1)^2}{3p-1}}\frac{(4 \pi )^{\frac{3-p}{3}}}{3^{\frac{4p}{3(3p-1)}}}\ncapa_p(1+\kst \varepsilon)\right]^{\frac{3p-1}{2p(p-1)}} \abs{\Omega}^{-\frac{2}{3(p-1)}} \nonumber\\
 &{}\leq \frac{3-p}{2}(4 \pi )^{-\frac{p-1}{2p}} \ncapa_p^{-\frac{3(p-1)}{2p(3-p)}} \bigg(\frac{3-p}{p-1}\bigg)^{\frac{p-1}{2p}}\frac{(1+ \kst \varepsilon)^{\frac{3p-1}{2p(p-1)}}}{ (1-\varepsilon)^{\frac{2}{3(p-1)}}}.\label{eq:upper_boundedness}
\end{align}

On the other hand, let $v\colon \set{ \abs{x} \geq R(1)} \subset \R^n \to (0,1]$ be the function such that $\set{v=t}= \set{\abs{x}=R(t)}$ and $\abs{\Omega_t}=4 \pi R(t)^3/3$.
Since by construction $\abs{\D v}= -4 \pi R(t)^{2}/V'(t)$, the function~$v$ is locally Lipschitz. By coarea formula, we have
\begin{align}
 \int_0^1 \frac{V(t)^{\frac{2p}{3}}}{(-V'(t))^{p-1}} \dif t &{}= \frac{1}{(36 \pi )^{\frac{p}{3}}}\int_0^1 { \int_{\set{v=t}} \abs{\D v}^{p-1} \dif \sigma} \dif t = \frac{1}{(36 \pi )^{\frac{p}{3}}} \int_{\set{ \abs{x}\geq R(1)}} \abs{\D v}^p\dif x \nonumber\\
 &{}\geq \frac{(4\pi)^{\frac{3-p}{3}}}{3^{\frac{2p}{3}}}\bigg(\frac{3-p}{p-1}\bigg)^{p-1} \ncapa_p(\set{\abs{x}=R(1)}) \nonumber\\
 &{}=\frac{1}{3^{p-1}}\bigg(\frac{3-p}{p-1}\bigg)^{p-1} \abs{\Omega}^{\frac{3-p}{3}}.\label{eq:lower_boundedness}
\end{align}
Plugging~\eqref{eq:lower_boundedness} and~\eqref{eq:upper_boundedness} into~\eqref{eq:integration_sharpcap}, we conclude the proof by arbitrariness of $\varepsilon$. \end{proof}

We are ready to prove the claimed upper bound of the $p$-isocapacitary mass in terms of the isoperimetric mass.

\begin{Theorem}\label{thm:upperbound_miso}
Let $(M,g)$ a $\CS^0$-asymptotically flat Riemannian $3$-manifold with possibly empty compact boundary. Then, for every $1<p\leq 2$, we have that
\[
 \ma^{\tp}_{\iso} \leq \ma_{\iso}.
\]
\end{Theorem}

\begin{proof} Let $(\Omega_j)_{j \in \N}$ be an exhaustion of $(M,g)$, then $\abs{\Omega_j}\to +\infty$ as $j \to +\infty$. In particular, by Theorem~\ref{thm:sharp-isopcapacitary}, we have that
\[
 \frac{1}{\ncapa_p(\partial \Omega_j)^{\frac{2-p}{3-p}}} \Bigg[\bigg(\frac{3 \abs{\Omega_j}}{4 \pi }\bigg)^{\frac{3-p}{3}} - \ncapa_p(\partial\Omega_j) \Bigg] \leq \frac{p(3-p)}{2}\ma(1+o(1))
\]
as $j\to +\infty$, where $\ma\in \R\cup\set{+\infty}$ is such that $\ma\geq \ma_{\iso}$. Hence
\[
 \limsup_{j \to +\infty}\frac{1}{\ncapa_p(\partial \Omega_j)^{\frac{2-p}{3-p}}} \Bigg[\bigg(\frac{3 \abs{\Omega_j}}{4 \pi }\bigg)^{\frac{3-p}{3}} - \ncapa_p(\partial\Omega_j) \Bigg] \leq \frac{p(3-p)}{2}\ma.
\]
Taking the supremum among all exhaustions $(\Omega_j)_{j \in \N}$, we conclude employing Proposition~\ref{prop:equivalence_of_iso_p_mass} for $\alpha = (3-p)/3$ and sending $\ma\to - \infty$ if $\ma_{\iso}=-\infty$. \end{proof}

Combining Lemma~\ref{cor:lowerbound_miso} and Theorem~\ref{thm:upperbound_miso}, we directly get the convergence of the $p$-iso\-ca\-pac\-i\-tary masses to the isoperimetric mass as $p \to 1^+$.

\begin{Corollary}
Let $(M,g)$ be a $\CS^0$-asymptotically flat Riemannian $3$-manifold with nonnegative scalar curvature and smooth, compact, minimal, possibly empty boundary. Assume that~$(M,g)$ satisfies~\eqref{eq:energy_hp} for any $1 < p < 1 + \delta$, for a fixed $0 < \delta < 2$. Then
\[
 \lim_{p\to 1^+} \ma^{\tp}_{\iso} = \ma_{\iso}.
\]
\end{Corollary}

We are ready to prove, in the stronger $\CS^1_\tau$-asymptotically flat assumptions, $\tau>1/2$, that the $p$-isocapacitary masses do actually coincide with each other.

\begin{proof}[Proof of Theorem~\ref{thm:equivalence-masses}] Assume that $\ma_{\ADM}$ is finite, otherwise~\cite[Theorem 4.13]{benatti_isoperimetricriemannianpenroseinequality_2022} yields ${\ma_{\iso}=+\infty}$ and Lemma~\ref{cor:lowerbound_miso} implies $\ma^{\tp}_{\iso}=+\infty$. The first inequality, under these assumptions, is the content of~\cite[Theorem 4.13]{benatti_isoperimetricriemannianpenroseinequality_2022}. Following the same lines of~\cite[Proposition 14]{jauregui_admmasscapacityvolumedeficit_2020} (based on computations contained in~\cite{fan_largespheresmallspherelimitsbrownyork_2009}, which in fact only relies on the $\CS^1$-character of the metric), we have
\begin{gather*}
 \frac{1}{16 \pi} \int_{\partial B_r} \H^2 \dif \sigma =1- \frac{2 \ma_{\ADM}}{r} + o\big(r^{-1}\big),\\
 \abs{ \partial B_r} = 4 \pi r^2 + 4 \pi \eta(r) + o(r),\\
 \frac{3 \abs{ B_r}}{4 \pi } = r^3 +\frac{3 \ma_{\ADM}}{2}r^2 + \frac{3 }{2} \eta(r) r + o\big(r^2\big) ,
\end{gather*}
as $r \to +\infty$, where $B_r= \set{ \abs{x} \leq r}$ and $\abs{\eta(r)}\leq \kst r^{2-\tau}$. Employing Proposition~\ref{thm:p-braymiao} and using Taylor's expansion of $\leftidx{_2}{F}{_1}$ around $0$ (see~\eqref{eq:hypergeometric_relation}), we have
\begin{align*}
\ncapa_p(\partial B_r) &{}\leq \big( r^2 + \eta(r) + o(r)\big)^{\frac{3-p}{2}}\leftidx{_2}{F}{_1}\bigg(\frac{1}{2}, \frac{3-p}{p-1}, \frac{2}{p-1};\frac{ 2\ma_{\ADM}}{r}+o\big(r^{-1}\big)\bigg)^{-(p-1)} \\
&{}\leq \big( r^2 + \eta(r) + o(r)\big)^{\frac{3-p}{2}}\bigg(1+ \frac{3-p}{2r} \ma_{\ADM} +o\big(r^{-1}\big) \bigg)^{-(p-1)} \\
&{} =r^{3-p}\bigg( 1+ \frac{3-p}{2r^2}\eta(r)+ o\big(r^{-1}\big)\bigg) \bigg( 1- \frac{(3-p)(p-1)}{2r}\eta(r)+ o\big(r^{-1}\big)\bigg) \\
&{} = r^{3-p}+ \frac{3-p}{2} \eta(r)r^{1-p}- \frac{(3-p)(p-1)}{2} r^{2-p} \ma_{\ADM} + o\big(r^{2-p}\big).
\end{align*}
Proposition~\ref{prop:equivalence_of_iso_p_mass} for $\alpha= (3-p)/3$ gives
\begin{align*}
 \ma^{\tp}_{\iso}& {}\geq \limsup_{r \to +\infty}\frac{2\ncapa_p(\partial B_r)^{\frac{p-2}{3-p}}}{p(3-p)}\Bigg( \bigg( \frac{3 \abs{ B_r}}{4 \pi }\bigg)^{\frac{3-p}{3}} - \ncapa_p(\partial B_r)\Bigg)\\
 & {}\geq \limsup_{r \to +\infty}\frac{2\ncapa_p(\partial B_r)^{\frac{p-2}{3-p}}}{p(3-p)} \bigg(\frac{p(3-p)}{2}r^{2-p}\ma_{\ADM}+ o (r^{2-p}) \bigg)= \ma_{\ADM},
\end{align*}
where the last identity is given by~\cite[Lemma 2.21]{benatti_asymptoticbehaviourcapacitarypotentials_2022}. The conclusion follows by Theorem~\ref{thm:upperbound_miso}, since $\ma_{\iso}= \ma_{\ADM}$~\cite[Theorem 4.13]{benatti_isoperimetricriemannianpenroseinequality_2022}. \end{proof}

\appendix

\section[Monotonicities along the p-capacitary potential]{Monotonicities along the $\boldsymbol{p}$-capacitary potential}
\label{app:monotonicities}
Here we slightly improve the monotonicity results in~\cite{agostiniani_riemannianpenroseinequalitynonlinear_2022,chan_monotonicitygreenfunctions_2022}. Inspired by these two works, we are approximating the $p$-capacitary potential with a family of smooth functions. To enter more in detail, let $(M,g)$ be a strongly $p$-nonparabolic Riemannian manifold with (possibly empty) boundary. Let $\Omega \subset M$ be homologous to $\partial M$ and $u_p$ is the solution to~\eqref{eq:pb-intro} starting at $\Omega$. For every $T>1$ let $\Omega_T$ be strictly homologous to $\partial M$ with connected boundary and containing $\set{u_p>\alpha_p(T)}$, where $\alpha_p(T) = T^{-(3-p)/(p-1)}$. Then, we define $u_p^\varepsilon$ as the solution to the following boundary value problem:
\begin{equation}\label{eq:pb-intro-eps}
 \begin{cases}
 \Delta_{p}^\varepsilon u^\varepsilon_p = 0 &\text{on $\Int\Omega_T \smallsetminus \Omega $,}\\
 u_p^\varepsilon=1 & \text{on $\partial \Omega$,}\\
 u_p^\varepsilon=u_p & \text{on $\partial\Omega_T$,}
 \end{cases}
\end{equation}
where
\[
 \Delta_{p}^\varepsilon f = \div\big( \abs{ \D f}_\varepsilon^{p-2} \D f\big) \qquad \text{and} \qquad {\abs{}}_\varepsilon= \sqrt{ \abs{}^2 +\varepsilon^2}.
\]
The function $u^\varepsilon_p$ is smooth away from the exterior boundary and converges in $\CS^{1,\beta}_{\loc}$ to the $p$-capacitary potential $u_p$ as $\varepsilon \to 0^+$. Indeed, this family was used in~\cite{dibenedetto_alphalocalregularityweak_1983} to prove $\CS^{1,\beta}_{\loc}$-regularity of $p$-harmonic functions. Moreover, looking more carefully at the proof of~\cite[Lemma~2.1]{lou_singularsetslocalsolutions_2008}, $\abs{\D u^\varepsilon_p}^{p-1}$~is uniformly bounded in $W^{1,2}_{\loc}$. Hence, up to a not relabeled subsequence, we can always assume that $\abs{ \D u_p^\varepsilon}^{p-1}$ weakly converges in $W^{1,2}_{\loc}$. Moreover, since $\abs{ \D u_p^\varepsilon}$ converges uniformly to~$| \D u_p |$, the weak limit of $\D \abs{ \D u_p^\varepsilon}^{p-1}$ must be $\D | \D u_p |^{p-1}$.

We are now ready to prove Theorem~\ref{thm:AMMO_monotonicity}.

\begin{proof}[Proof of Theorem~\ref{thm:AMMO_monotonicity}] For ease of computations, we rewrite the function $t\mapsto \ma_{H}^{\tp}(\partial \Omega_t)$ in terms of the $p$-capacitary potential, that is,
\[
 U_p(t)= 4 \pi t + \frac{(p-1)^2}{(3-p)^2} t^{\frac{5-p}{p-1}} \int_{\set{u_p = \alpha_p(t)}} | \D u_p |^2 \dif \sigma -\frac{(p-1)}{(3-p)} t^{\frac{2}{p-1}} \int_{\set{u_p = \alpha_p(t)}} | \D u_p |\H \dif \sigma.
\]
Observe that
\[
 \frac{\ncapa_p(\partial \Omega)^\frac{1}{3-p}}{8\pi} U_p\big(\ee^{\frac{t}{3-p}}\big) = \ma_{H}^{\tp}(\partial \Omega_t)= \frac{1}{8\pi}F_p\big(\ncapa_p(\partial \Omega)^\frac{1}{3-p} \ee^{\frac{t}{3-p}}\big),
\]
where $F_p$ is the function defined in~\cite{agostiniani_riemannianpenroseinequalitynonlinear_2022}. Equation~\eqref{eq:AMMO_monotonicity} now follows from computations in~\cite[Section 1.2]{agostiniani_riemannianpenroseinequalitynonlinear_2022}. The function $U_p$ is well defined on $[0,+\infty)$. Indeed, one can observe that
\[
 \abs{ \H} | \D u_p |^2 = | \D u_p |^{3-p} | \D | \D u_p |^{p-1}| \in L^2_{\loc}(M\smallsetminus \Omega),
\]
since $1<p<3$ and $| \D u_p |^{p-1} \in L^\infty_{\loc}\cap W^{1,2}_{\loc} (M\smallsetminus \Omega)$ by~\cite[Lemma~2.1]{lou_singularsetslocalsolutions_2008}. Then, by coarea formula, the function
\[
 t \mapsto \int_{\set{u_p= \alpha_p(t)}} | \D u_p | \H \dif \sigma\in L^1_{\loc}(0,+\infty)
\]
and its equivalence class does not depend on the representative of $| \D u_p | \H $. In particular, the function $t\mapsto \ma_{H}^{\tp}(\partial \Omega_t) \in L^1_{\loc}(0,+\infty)$.

It only remains to prove that it has nonnegative first derivative in the sense of distributions, which both gives that $t\mapsto \ma_{H}^{\tp}(\partial \Omega_t) \in \BV_{\loc}(0,+\infty)$ and that admits a nondecreasing representative. Fix $T>0$ and $u^\varepsilon_p$ be the solution to~\eqref{eq:pb-intro-eps}. One can now define the function
\[
 F_p^\varepsilon(t)= 4 \pi t + \frac{(p-1)^2}{(3-p)^2} t^{\frac{5-p}{p-1}} \int_{\set{u^\varepsilon_p = \alpha_p(t)}} \abs{ \D u^\varepsilon_p}^2 \dif \sigma -\frac{(p-1)}{(3-p)} t^{\frac{2}{p-1}} \int_{\set{u^\varepsilon_p = \alpha_p(t)}} \abs{ \D u^\varepsilon_p}\H^\varepsilon\dif \sigma
\]
for every $t \in (1,T)$, where $\H^\varepsilon$ is the mean curvature of the level $\set{u^\varepsilon_p = \alpha_p(t)}$. By~\cite[Lemma~1.2]{agostiniani_riemannianpenroseinequalitynonlinear_2022}, we have that the function $F^\varepsilon_p$ is almost monotone, in the sense that
\[
 F^\varepsilon_p(t)- F^\varepsilon_p(s) \geq -\varepsilon
 \frac{(p+1)^2(p-1)}{(3-p)^3} \int_{ \set{\alpha_p(t) < u^\varepsilon_p<\alpha_p(s)}}\abs{ \D u_p^\varepsilon}^2(u^\varepsilon_p)^{-\frac{p-1}{3-p}-3} \dif \mu
\]
holds for almost every $0\leq s\leq t \leq T$.
In particular, by coarea formula and Lebesgue differentiation theorem,
\[
 (F^\varepsilon_p)'(t) \geq
 - \varepsilon \frac{(p+1)^2}{(3-p)^2}t^{\frac{2(3-p)}{p-1}} \int_{\set{u^\varepsilon_p=\alpha_p(t)}} \abs{ \D u^\varepsilon_p} \dif \sigma
\]
holds for almost every $t \in (1,T)$.

The monotonicity follows if one can prove the following claim.
\begin{Claim}
$F_p^\varepsilon$ converges to $F_p$ in the sense of distributions as $\varepsilon\to 0^+$.
\end{Claim}
Indeed, if that is the case, for every nonnegative test function $\varphi \in \CS^\infty_c(1,T)$, we obtain that
\begin{align*}
 -\int_1^T \varphi'(t) F_p(t)\dif t & =-\lim_{\varepsilon\to 0} \int_1^T\varphi'(t) F^p_\varepsilon(t)\dif t \\
 & \geq - \frac{(p+1)^2}{(3-p)^2} \lim_{\varepsilon \to 0} \varepsilon \int_1^T \varphi(t) t^{\frac{2(3-p)}{p-1}} \int_{\set{u^\varepsilon_p=\alpha_p(t)}} \abs{ \D u^\varepsilon_p} \dif \sigma\dif t \\
 & =-\frac{(p+1)^2}{(3-p)^2} \lim_{\varepsilon \to 0} \varepsilon \int_{\Omega_T \smallsetminus \Omega} \varphi \big(\alpha_p^{-1}\big(u^\varepsilon_p\big)\big) \frac{ \big| \D u^\varepsilon_p \big|^2}{\big(u^\varepsilon_p\big)^{\frac{3-p}{p-1}+3}} \dif \mu=0,
\end{align*}
since $u_p^\varepsilon$ converges to $u_p$ in $\CS^{1,\beta}_{\loc}$. This shows that $F_p$ has nonnegative first derivative in the sense of distributions, proving its monotonicity.

We now turn to prove the claim. Consider any $\varphi \in \CS^\infty_c(0,+\infty)$. The first term is independent of $\varepsilon$. As far as the second term is concerned, by coarea formula, we have that
\[
 \int_1^T\varphi(t)t^{\frac{5-p}{p-1}} \int_{\set{u_p^\varepsilon=\alpha_p(t)}} \abs{\D u^\varepsilon_p}^2 \dif \sigma \dif t = \frac{(p-1)}{(3-p)} \int_{\Omega_T \smallsetminus \Omega} \big(u^\varepsilon_p\big)^{-\frac{7-p}{3-p}}\varphi \big(\alpha_p^{-1}\big(u^\varepsilon_p\big)\big) \abs{ \D u_p^\varepsilon}^3 \dif \mu.
\]
Since the function $\varphi$ is smooth with compact support in $(1,T)$ and $u_p^\varepsilon$ converges to $u_p$ in $\CS^{1,\beta}_{\loc}$, the right-hand side converges to
\[
 \frac{(p-1)}{(3-p)} \int_{\Omega_T \smallsetminus \Omega} u_p^{-\frac{7-p}{3-p}}\varphi \big(\alpha_p^{-1}\big(u^\varepsilon_p\big)\big)| \D u_p |^3 \dif \mu = \int_1^T \varphi(t)t^{\frac{5-p}{p-1}} \int_{\set{u_p^\varepsilon=\alpha_p(t)}} |\D u_p |^2 \dif \sigma \dif t,
\]
where the identity follows by coarea formula. The last term is a little trickier since it involves second derivatives of the function $u_p^\varepsilon$ that are not converging uniformly as $\varepsilon \to 0$ to the corresponding ones for $u_p$. Employing again the coarea formula and by straightforward computations, we then have that
\begin{gather*}
 \int_1^T\varphi(t)t^{\frac{2}{p-1}} \int_{\set{u_p^\varepsilon=\alpha_p(t)}} \abs{\D u^\varepsilon_p} \H_\varepsilon\dif \sigma \dif t \\
 \qquad {}=
 \int_1^T\varphi(t)t^{\frac{2}{p-1}} \int_{\set{u_p^\varepsilon=\alpha_p(t)}} \frac{\big\langle \D \big| \D u^\varepsilon_p\big|^{p-1}\big|\D u^\varepsilon_p\big\rangle}{\big| \D u_p^\varepsilon\big|^{p-1} }\Bigg( 1 - \frac{(p-2)}{(p-1)} \frac{ \varepsilon^2}{\big|\D u^\varepsilon_p\big|_\varepsilon^2}\Bigg)\dif \sigma \dif t\\
 \qquad {}=\frac{(p-1)}{(3-p)} \int_{\Omega_T \smallsetminus \Omega} \big(u^\varepsilon_p\big)^{-\frac{4}{3-p}}\varphi \big(\alpha_p^{-1}\big(u^\varepsilon_p\big)\big)\frac{\big\langle\D \big| \D u^\varepsilon_p\big|^{p-1}\big|\D u^\varepsilon_p\big\rangle}{\big| \D u_p^\varepsilon\big|^{p-2} }\Bigg( 1 - \frac{(p-2)}{(p-1)} \frac{ \varepsilon^2}{\big|\D u^\varepsilon_p\big|_\varepsilon^2}\Bigg) \dif \mu.
\end{gather*}
As before we want to prove that the right-hand side converges to the corresponding term for $u_p$. Since $u^\varepsilon_p\to u_p$ in $\CS^{1, \beta}_{\loc}$ and $\D\big| \D u^\varepsilon_p\big|^{p-1}\rightharpoonup \D| \D u_p |^{p-1} $ weakly in $L^2_{\loc}$, we have that
\begin{gather*}
 \lim_{\varepsilon \to 0}\frac{(p-1)}{(3-p)} \int_{\Omega_T \smallsetminus \Omega} \big(u^\varepsilon_p\big)^{-\frac{4}{3-p}}\varphi \big(\alpha_p\big(u^\varepsilon_p\big)\big) \frac{\big\langle\D \big| \D u^\varepsilon_p\big|^{p-1}\big|\D u^\varepsilon_p\big\rangle}{\big| \D u_p^\varepsilon\big|^{p-2} }\dif \mu\\
\qquad {} =\frac{(p-1)}{(3-p)} \int_{\Omega_T \smallsetminus \Omega} (u_p)^{-\frac{4}{3-p}}\varphi \big(\alpha_p^{-1}\big(u^\varepsilon_p\big)\big)\frac{\big\langle\D | \D u_p |^{p-1}|\D u_p\big\rangle}{| \D u_p |^{p-2} }\dif \mu\\
\qquad {} =\int_1^T\varphi(t)t^{\frac{2}{p-1}}\int_{\set{u_p=\alpha_p(t)}} |\D u_p | \H \dif \sigma \dif t.
\end{gather*}
Moreover, the remaining term vanishes. H\"older's inequality and equi-boundedness in $L_{\loc}^2(\Omega_T \smallsetminus \Omega)$ of $\big|\D \big|\D u^\varepsilon_p\big|^{p-1}\big|$ yield
\begin{equation}\label{eq:main_inequality_monotonicity_AMMO}
 \int_{K}\abs{\D u^\varepsilon_p}^{3-p}\big|\D \abs{ \D u^\varepsilon_p}^{p-1}\big|\frac{ \varepsilon^2}{\big| \D u^\varepsilon_p\big|_\varepsilon^2}\dif \mu\leq \kst_1\Bigg(\int_K\frac{\varepsilon^4}{\big| \D u^\varepsilon_p\big|_\varepsilon^4} \abs{\D u^\varepsilon_p}^{6-2p}\dif \mu\Bigg)^\frac{1}{2}
\end{equation}
for every $K$ compactly contained in $\set{ \alpha_p(T) < u_p^\varepsilon <1}$ and for some positive constant $\kst_1$. Observe that
\[
 \abs{\D u^\varepsilon_p}^{6-2p}\frac{\varepsilon^4}{\big| \D u^\varepsilon_p\big|_\varepsilon^4} \leq \abs{\D u^\varepsilon_p}^{6-2p}\leq \kst_2,
\]
since the function $\big| \D u_\varepsilon^p\big|$ converges locally uniformly and $1<p<3$. The left-hand side converges almost everywhere to $0$. Indeed, if a point belongs to critical set of $u_p$, \smash{$\big|\D u_\varepsilon^p\big|^{6-2p}\to 0$} as $\varepsilon\to 0$. Otherwise, $\big| \D u_\varepsilon^p\big|$ is definitely bounded away from $0$, then \smash{$\big| \D u^\varepsilon_p\big|_\varepsilon^4$} is not vanishing, thus the left-hand side is controlled by $\varepsilon^4$ up to a constant. By dominated convergence theorem, the right-hand side in~\eqref{eq:main_inequality_monotonicity_AMMO} approaches $0$ as $\varepsilon \to 0$, so does the left-hand side, concluding the step. \end{proof}

We use this theorem to study the quantity introduced in~\cite[Theorem 1.2]{chan_monotonicitygreenfunctions_2022} along the level set of the solution $w_p$ to~\eqref{eq:pb-intro}.

\begin{Theorem}\label{thm:CCLT_monotonicity}
 Let $(M,g)$ be a strongly $p$-nonparabolic Riemannian $3$-manifold with smooth, compact and connected possibly empty boundary $\partial M$. Assume that $H_2(M, \partial M;\Z) = \set{0}$. Let $\Omega\subseteq M$ be bounded closed with connected $\CS^1$-boundary homologous to $\partial M$ and with $\h\in L^2(\partial \Omega)$. Let $w_p$ be the solution to~\eqref{eq:pb-intro} starting at $\Omega$. Then, denoting $\Omega_t = \set{w_p \leq t}$, the function
 \begin{equation}\label{eq:monotoncity_chan}
 t \mapsto \frac{\ncapa_p(\partial \Omega_t)^{-\frac{1}{p-1}}}{4\pi(3-p)}\Bigg( 4 \pi - \int_{\partial \Omega_t} \frac{\abs{\D w_p}^2}{(3-p)^2}\dif \mu\Bigg)
 \end{equation}
 belongs to $W^{1,1}_{\loc}(0,+\infty)$ and
 \begin{equation}\label{eq:weak_derivative}
 \frac{\dif}{\dif t}\Bigg[\ncapa_p(\partial \Omega_t)^{-\frac{1}{p-1}}\Bigg( 4 \pi - \int_{\partial \Omega_t} \frac{\abs{\D w_p}^2}{(3-p)^2}\dif \sigma\Bigg)\Bigg] =- \frac{8\pi}{p-1}\ncapa_p(\partial \Omega_t)^{-\frac{2}{(3-p)(p-1)}}\ma_{H}^{\tp}(\partial \Omega_t),
 \end{equation}
 for almost every $t \in [0,+\infty)$.
\end{Theorem}

\begin{Remark}\label{rmk:cclt_mon}
 Observe that~\eqref{eq:monotoncity_chan} is, up to multiplying by a constant and changing variables, the quantity studied in~\cite{chan_monotonicitygreenfunctions_2022,munteanu_comparisontheorems3dmanifolds_2023} for $p\in(1,2]$ along the level set of the $p$-Green function. In particular, coupled with Theorem~\ref{thm:AMMO_monotonicity}, the above result implies that the monotonicity property of~\eqref{eq:monotoncity_chan} is preserved if, in place of the $p$-Green's function, the $p$-capacitary potential of a connected $\partial \Omega$ with nonnegative $p$-Hawking mass is considered. Clearly, the monotonicity results in~\cite{chan_monotonicitygreenfunctions_2022, munteanu_comparisontheorems3dmanifolds_2023} are recovered applying Theorem~\ref{thm:CCLT_monotonicity} to the $p$-Green's function. Hence, we settle the question raised in~\cite{chan_monotonicitygreenfunctions_2022} about the monotonicity of their quantities in the range $p \in (2, 3)$.

 Finally, observe that $\tilde{\ma}^{\tp}_H$ is obtained multiplying the function~\eqref{eq:monotoncity_chan} by the $p$-capacity term $\ncapa_p(\partial \Omega_t)^{2/(3-p)(p-1)}$ which is exponentially growing as $t \to +\infty$. This term forces the quantity~$\tilde{\ma}^{\tp}_H$ to be monotone \emph{nondecreasing} under assumption~\eqref{eq:energy_hp} (see Lemma~\ref{thm:inequality_modifiedgeroch_AMMO}) even when~\eqref{eq:monotoncity_chan} is monotone \emph{nonincreasing}.
\end{Remark}

\begin{proof}[Proof of Theorem~\ref{thm:CCLT_monotonicity}] Fix $T>1$ and let $u^\varepsilon_p$ the solution to the problem~\eqref{eq:pb-intro-eps}. Consider any $\varphi \in \CS^\infty_c(0,T)$. Employing coarea formula and integration by parts, we have that
\begin{align}
 \int_1^{T} \varphi'(t) \int_{\set{u^\varepsilon_p=\alpha_p(t)}} \abs{ \D u^\varepsilon_p}^2 \dif \sigma \dif t &{}= \frac{(p-1)}{(3-p)} \int_{\Omega_T \smallsetminus \Omega} \varphi' \big(\big(u_p^\varepsilon\big)^{-\frac{p-1}{3-p}} \big)\big(u_p^\varepsilon\big)^{-\frac{2}{3-p}} \abs{ \D u^\varepsilon_p}^3 \dif \mu \nonumber\\
 &{}=- \int_{\Omega_T \smallsetminus \Omega} \big\langle\big| \D u^\varepsilon_p\big|\D u^\varepsilon_p\big| \D \big[\varphi \big(\big(u_p^\varepsilon\big)^{-\frac{p-1}{3-p}}\big)\big]\big\rangle \dif \mu \nonumber\\
 &{}= \int_{\Omega_T \smallsetminus \Omega} \div \big( \abs{ \D u^\varepsilon_p} \D u^\varepsilon_p\big) \varphi \big(\big(u_p^\varepsilon\big)^{-\frac{p-1}{3-p}}\big) \dif \mu. \label{eq:main_chan}
\end{align}
Clearly, employing $\CS^{1, \beta}_{\loc}$ convergence of $u^\varepsilon_p \to u_p$ as $\varepsilon \to 0$, we obtain that
\[
 \lim_{\varepsilon\to 0}\int_1^{T} \varphi'(t) \int_{\set{u^\varepsilon_p=\alpha_p(t)}} \abs{ \D u^\varepsilon_p}^2 \dif \sigma \dif t = \int_1^{T} \varphi'(t) \int_{\set{u_p=\alpha_p(t)}} | \D u_p |^2 \dif \sigma \dif t.
\]
Moreover, a straightforward computation leads to
\[
 \div \big( \abs{ \D u^\varepsilon_p} \D u^\varepsilon_p\big) =\frac{\big| \D u^\varepsilon_p\big|^{2-p}}{(p-1)} \Bigg( (3-p) \frac{ \big| \D u^\varepsilon_p\big|^2}{\big| \D u^\varepsilon_p \big|^2_\varepsilon}- \frac{ \varepsilon^2}{\big| \D u^\varepsilon_p\big|_{\varepsilon}^2} \Bigg) \big\langle\D \big| \D u^\varepsilon_p \big|^{p-1} \big| \D u^\varepsilon_p\big\rangle.
\]
Arguing as in the previous theorem, since $\big|\D u^\varepsilon_p\big| \to |\D u_p |$ locally uniformly and $\D\big|\D u^\varepsilon_p\big|^{p-1} \rightharpoonup\D |\D u_p |^{p-1}$ weakly $L^2_{\loc}$ as $\varepsilon\to 0$, one gets
\begin{gather*}
 \lim_{\varepsilon \to 0} \int_{\Omega_T \smallsetminus \Omega} \frac{\big| \D u^\varepsilon_p\big|^{2-p}}{(p-1)} \Bigg( (3-p) \frac{ \big| \D u^\varepsilon_p\big|^2}{\big| \D u^\varepsilon_p \big|^2_\varepsilon}- \frac{ \varepsilon^2}{\big| \D u^\varepsilon_p\big|_{\varepsilon}^2} \Bigg) \big\langle\D \big| \D u^\varepsilon_p \big|^{p-1} \big| \D u^\varepsilon_p\big\rangle \varphi \big(\big(u_p^\varepsilon\big)^{-\frac{p-1}{3-p}}\big) \dif \mu\\
 \qquad {}=
 \frac{(3-p)}{(p-1)} \int_{\Omega_T \smallsetminus \Omega} \frac{ \big\langle \D | \D u_p |^{p-1} | \D u_p\big\rangle}{| \D u_p |^{p-2}} \varphi \big((u_p)^{-\frac{p-1}{3-p}}\big)\dif \mu \\
 \qquad {}=
 \frac{(3-p)^2}{(p-1)^2}\int_1^{T}\varphi(t) t^{-\frac{2}{p-1}} \int_{\Omega_T \smallsetminus \Omega} \H \abs{ \D u} \dif \sigma \dif t,
\end{gather*}
where the last equality follows by coarea formula and $\H= \big\langle\D | \D u_p |^{p-1} | \D u_p\big\rangle / | \D u_p |^p$. Then, passing to the limit as $\varepsilon \to 0$ in~\eqref{eq:main_chan}, we obtain
\[
 \int_1^{T} \varphi'(t) \int_{\set{u_p=\alpha_p(t)}} | \D u_p |^2 \dif \sigma \dif t =
 \frac{(3-p)^2}{(p-1)^2}\int_1^{T}\varphi(t) t^{-\frac{2}{p-1} } \int_{\set{u_p= \alpha_p(t)}} \H \abs{ \D u} \dif \sigma \dif t.
\]
By arbitrariness of $T$ and $\varphi$ one has that the function
\[
 t \mapsto H_p(t)=\Bigg(4 \pi t^{-\frac{3-p}{p-1}}-\frac{(p-1)^2}{(3-p)^2}t^{\frac{3-p}{p-1}} \int_{\set{u=\alpha_p(t)}} |\D u_p |^2 \dif \sigma\Bigg)
\]
belongs to $W^{1,1}_{\loc}(1,+\infty)$ and its derivative is given by
\begin{align*}
H'_p(t)={}& -\frac{(3-p)}{(p-1)}t^{-\frac{p+1}{p-1}}\Bigg(4 \pi t+ \frac{(p-1)^2}{(3-p)^2} t^{\frac{5-p}{p-1}} \\
&{}\times \int_{\set{u=\alpha_p(t)}} |\D u_p |^2 \dif \sigma - \frac{(p-1)}{(3-p)} t^{\frac{2}{p-1}}\int_{\set{u_p= \alpha_p(t)}} \H \abs{ \D u} \dif \sigma\Bigg)
\end{align*}
holds for almost every $t \in (1,+\infty)$. Then
\begin{align*}
 \frac{\dif}{\dif t}\Bigg[\ncapa_p(\partial \Omega_t)^{-\frac{1}{p-1}}\Bigg( 4 \pi - \int_{\partial \Omega_t} \frac{\big|\D w_p\big|^2}{(3-p)^2}\dif \mu\Bigg)\Bigg] &{}= \frac{\ncapa_p(\partial \Omega)^{-\frac{1}{p-1}}}{3-p} H'_p\big(\ee^{\frac{t}{3-p}}\big) \ee^{\frac{t}{3-p}}\\
 &{}=- \frac{8\pi}{p-1}\ncapa_p(\partial \Omega_t)^{-\frac{2}{(3-p)(p-1)}}\ma_{H}^{\tp}(\partial \Omega_t),
\end{align*}
concluding the proof.\end{proof}

\subsection*{Acknowledgements}
Part of this work has been carried out during the authors' attendance to the \emph{Thematic Program on Nonsmooth Riemannian and Lorentzian Geometry} that took place at the Fields Institute in Toronto. The authors warmly thank the staff, the organizers and the colleagues for the wonderful atmosphere and the excellent working conditions set up there.
L.B.~is supported by the European Research Council's (ERC) project n.853404 ERC VaReg -- \textit{Variational approach to the regularity of the free boundaries}, financed by the program Horizon 2020, by PRA\_2022\_11 and by PRA\_2022\_14.
M.F.~has been supported by the European Union -- NextGenerationEU and by the University of Padova under the 2021 STARS Grants@Unipd programme ``QuASAR''.
The authors are members of Gruppo Nazionale per l'Analisi Matematica, la Probabilit\`a e le loro Applicazioni (GNAMPA), which is part of the Istituto
Nazionale di Alta Matematica (INdAM), and are partially funded by the GNAMPA project ``Problemi al bordo e applicazioni geometriche''.
The authors are grateful to S.~Hirsch and F.~Oronzio for their interest in the work and for pleasureful and useful conversations on the subject. The authors warmly thank the anonymous referees for their thorough reading of the paper, and for the precious suggestions that allowed to improve the quality of the paper.

\pdfbookmark[1]{References}{ref}
\LastPageEnding

\end{document}